\documentclass[11pt,twoside]{article}

\usepackage[english]{babel}
\usepackage{amsmath}
\usepackage{amsfonts}
\usepackage{amssymb}
\usepackage[all]{xy}
\usepackage{a4}
\usepackage{graphicx}
\usepackage{theorem}

\DeclareRobustCommand{\INSERT}[1]{%
   \ifmmode
      \mathbf{INSERT #1}%
   \else
      \textbf{INSERT #1}%
   \fi
}

\newcommand{\RR}{\mathbb{R}}

\newcommand{\ZZ}{\mathbb{Z}}
\newcommand{\BB}{\mathbb{B}}
\newcommand{\PP}{\mathbb{P}}
\newcommand{\EE}{\mathbb{E}}
\newcommand{\GG}{\mathbb{G}}
\newcommand{\KE}{\mathcal{E}}
\renewcommand{\gg}{\mathfrak{g}}
\newcommand{\FF}{\mathcal{F}}

\makeatletter
\renewcommand{\subsection}{\@startsection
{subsection} {1}
{.2cm}                   
{\baselineskip} {-\fontdimen2\font plus -\fontdimen3\font minus -\fontdimen4\font}
{\normalfont\normalsize\bfseries}} \makeatother

\newenvironment{proof}{\removelastskip\par\medskip
\noindent{\em Proof.} \rm}{\penalty-20\null\hfill$\square$\par\medbreak}

\begin{document}

\title{On the classification of Regular Lie groupoids}
\author{I. Moerdijk \\ University of Utrecht}
\date{}

\maketitle

\begin{abstract}
We observe that any regular Lie groupoid $G$ over an manifold $M$ fits into an extension $K \rightarrow G
\rightarrow E$ of a foliation groupoid $E$ by a bundle of connected Lie groups $K$. If $\FF$ is the foliation on
$M$ given by the orbits of $E$ and $T$ is a complete transversal to $\FF$, this extension restricts to $T$, as an
extension $K_{T}\rightarrow G_{T}\rightarrow E_{T}$ of an \'etale groupoid $E_{T}$ by a bundle of connected groups
$K_{T}$. We break up the classification into two parts. On the one hand, we classify the latter extensions of
\'etale groupoids by (non-abelian) cohomology classes in a new \v{C}ech cohomology of \'{e}tale groupoids. On the
other hand, given $K$ and $E$ and an extension $K_{T}\rightarrow G_{T}\rightarrow E_{T}$ over $T$, we present a
cohomological obstruction to the problem of whether this is the restriction of an extension $K \rightarrow G
\rightarrow E$ over $M$; if this obstruction vanishes, all extensions $K \rightarrow G \rightarrow E$ over $M$
which restrict to a given extension over the transversal together form a principal bundle over a ``group'' of
bitorsors under $K$.

\vspace{.5cm}
\end{abstract}
\newpage
\tableofcontents
\newpage

\section*{Introduction}
A groupoid is a category in which every arrow is invertible; a Lie groupoid is such a groupoid with an additional
smooth structure. The notion of Lie groupoid at the same time generalizes that of smooth manifolds and of Lie
groups, and goes back as far as to the work of Ehresmann in the fifties. Lie groupoids now play a central r\^ole
in several areas of mathematics, for example in the theory of foliations and non-commutative geometry, and in
quantization and deformation theory; see e.g. the books by Connes, Canas da Silva-Weinstein, Mackenzie(to appear)
and Landsman.

A Lie groupoid is given by two manifolds, a manifold $M$ (of ``objects'') and a manifold $G$ (of ``arrows'')
together with several smooth structure maps, including the source and the target map $\xymatrix{s,t: G
\ar@<.5ex>[r] \ar@<-.5ex>[r] & M}$. For two points $x$ and $y$ in $M$, an arrow from $x$ to $y$ in the groupoid is
simply a point $g\in G$ with $s(g)=x$ and $t(g)=y$. For a given $x$, the set $t s^{-1}(x)$ of all targets of
arrows out of $x$ is an immersed submanifold of $M$, called the \textit{orbit} of $G$ through $x$. The set
$G_{x}\subset G$ of all arrows from $x$ to itself is a Lie group, called the isotropy group of $G$ at $x$.

It seems impossible to classify in some systematic and informative way all possible Lie groupoids over a given
manifold $M$ of objects. In this paper, we will give a classification of those groupoids $G$ over $M$ which are
\textit{regular}. This means that the orbits of $G$ (or equivalently, its isotropy groups) have a constant
dimension. Although this is a significant restriction on a Lie groupoid, many important classes of Lie groupoids
are regular. For example, groupoids arising in foliation theory are always regular, as are the Lie groupoids
involved in integration and quantization of regular Poisson manifolds \cite{AlcHec}. Much of the work of Mackenzie
is concerned with the study and application of transitive Lie groupoids, and these are all regular as well.
Furthermore, even when a groupoid $G$  over $M$ is not regular, it is often possible to provide a stratification
of $M$ for which the groupoid restricts to a regular groupoid on each of the strata (see e.g. \cite{Nistor}).

An important step in the understanding of regular groupoids has been made in a recent paper by Weinstein \cite{W}.
For those regular groupoids which are proper (these have compact isotropy groups), Weinstein provides a normal
form locally around each orbit.

In this paper, we take a viewpoint complementary to that of Weinstein, and adopt as a starting point the
restriction of a given regular groupoid to submanifolds of $M$ which are transversal to the orbits of $M$.

More explicitly, for a given regular groupoid $G$ over $M$, the orbits of $G$ define a foliation $\FF$ of $M$. Let
$T\subset M$ be a submanifold of $M$ which is transversal to this foliation and meets every leaf; in other words,
$T$ is a ``complete'' transversal. For such a transversal $T$, the groupoid $G$ restricts to a Morita equivalent
groupoid $G_{T}$ over $T$. Our first aim will be to study the problem to which extent $G$ is determined by
$G_{T}$. Our solution to this problem takes the following form. First, the groupoid $G$ fits into an extension
\begin{equation}\label{rr1}
K\rightarrow G\rightarrow E,
\end{equation}
where $K$ is a bundle of connected Lie groups over $M$, and $E$ is a foliation groupoid over $M$ (this is a
groupoid which ``integrates'' $\FF$, such as the holonomy groupoid of $\FF$). Similarly, $G_{T}$ fits into an
extension
\begin{equation}\label{rr2}
K_{T}\rightarrow G_{T}\rightarrow E_{T}
\end{equation}
over the transversal $T$. Not every such extension (\ref{rr2}) comes from an extension (\ref{rr1}) over $M$. In
fact, we show that there is an obstruction in $H^{2}(M,\underline{K})$, where $\underline{K}$ is the sheaf of
sections of the bundle $K$ (cf. Proposition 5.1 and Proposition 10.3). If this obstruction vanishes, there is a
beautiful description of all possible extension (\ref{rr1}) which restrict to the same extension (\ref{rr2}) over
$T$: they form a principal bundle under a group of bitorsors over $M$ (Theorem 6.6). In the special case that the
groupoid $G$ is transitive, the transversal $T$ reduces to a single point, and these results specialize to a
variation of the results in \cite{Mack3}.

Our next aim is to classify extensions of the form (\ref{rr2}) by degree 2 cohomology classes, much as in the
classical case of extensions of (discrete) groups, going back to \cite{MacLane}, \cite{Dedecker}, and others. At
first sight this would seem impossible, because of the lack of cross sections of the quotient map
$G_{T}\rightarrow E_{T}$ in (\ref{rr2}). However, the extensions of the form (\ref{rr2}) have a special property
implied by transversality of $T$, namely that the groupoid $E_{T}$ is \'etale. Etale groupoids have played a
central r\^ole in foliation theory ever since the early work of Haefliger (\cite{H2}), and are much better
understood than general Lie groupoids. In particular, any \'etale groupoid admits a description of the cohomology
of its classifying space in terms of the cohomology of a small (discrete) category, and it is this so-called
``embedding'' category which can be used to circumvent the problem of lack of cross sections. We will define a
relative (non-abelian) \v{C}ech cohomology of \'etale groupoids by using these small embedding categories (in \S7
and \S11), and prove that extensions of the form (\ref{rr2}) correspond to degree 2 cohomology classes (see
Theorems 8.2 and 12.3). This classification is closely related to the theory of gerbes on manifolds and on \'etale
groupoids, and it relies on Giraud's non-abelian cohomology \cite{Giraud} and subsequent work of Breen
\cite{Breen}. Gerbes over \'etale groupoids occur in work of Lupercio and Uribe \cite{LupUrb}.

To conclude this introduction, we make some remarks about our order of presentation. This paper is divided into
four parts. The first part provides some preliminaries on regular Lie groupoids. It serves to fix notation and
terminology, to recall some basic facts, and to motivate the subsequent results. The second part is concerned with
the problem of how to enlarge extensions of the form (\ref{rr2}) over $T$ to ones of the form (\ref{rr1}) over all
of $M$. This part is really a contribution to the general theory of extensions of Lie groupoids, and does not
depend on the special way the extensions (\ref{rr1}) and (\ref{rr2}) arise. (For example the fact that $E_{T}$ is
\'etale is never used here.) The third part introduces embedding categories for \'etale groupoids and the relative
``\v{C}ech'' cohomology necessary to obtain a cohomological classification of extensions of the form (\ref{rr2}).

In the presentation of this cohomological classification in the third part, as well as that of the obstruction
class in the second part, we initially restrict our attention to the case where the bundles of Lie groups involved
($K$ in (\ref{rr1}) and $K_{T}$ in (\ref{rr2})) are abelian. The general non-abelian case involves Giraud's theory
of gerbes and bands (``liens''), with which some readers interested in Lie groupoids may be less familiar. For
this reason, we have postponed the discussion of non-abelian cohomology to a separate fourth part. In this part,
we provide the necessary preliminaries for non-abelian cohomology of manifolds and of small categories,
respectively, in \S9 and \S11, and we show how to adapt the results of \S5 and \S8 to the non-abelian case in \S10
and \S12.

\textit{Acknowledgements:} I would like to thank L. Breen, A. Kumjian, M. Crainic, E. Lupercio, K. Mackenzie, J.
Mrcun and A. Weinstein for many helpful conversations.

The main results of this paper were first presented at the LMS Workshop on Quantization, etc. in Manchester (July
2001), and I thank M. Prest and A. Voronov for inviting me to speak there.

\newpage
\part{Regular Lie groupoids}

\section{Lie groupoids}
In this preliminary section we recall the basic definitions and fix some notation. For more background and many
examples, see \cite{Mack, Connes, Weinstein1}.

\subsection{Lie groupoids.}
Let $M$ be a smooth (Hausdorff) manifold. A \textit{Lie groupoid over} $M$  consists of a manifold  $G$ equipped
with several smooth structure maps. To begin with, there are two submersions
\[\xymatrix{s,t: G \ar@<.5ex>[r] \ar@<-.5ex>[r] & M} \]
called the source and the target maps. A point $g\in G$ with $s(g)=x$ and $t(g)=y$ will be referred to as an
\textit{arrow} from $x$ to $y$, and denoted by $g:x \rightarrow y$. Next, there is a smooth associative
composition operation $m: G \times_{M} G \rightarrow G$ mapping a pair ${g,h}$ with $s(g)=t(h)$ to the composition
$m(g,h)$, also simply denoted by $g h$ or $g \cdot h$. Finally, there is an involution $i: G \rightarrow G$
mapping each $g: x \rightarrow y$ in $G$ to a 2-sided inverse $i(g)=g^{-1}$ for composition, and a map $u:M
\rightarrow G$ mapping each $x\in M$ to a 2-sided unit $u(x)=1_{x}$ for composition. Thus, all in all, a Lie
groupoid is a system of manifolds and maps
\[(M,G,s,t,m,i,u),\]
for which we simply write $G$. As will be clear from the examples, we will allow the manifold $G$ to be
non-Hausdorff.

\subsection{Lie algebroids.}
The Lie algebroid $\gg$ of a Lie groupoid $G$ is the vector bundle $u^{*}(T^{s}(G))$ over $M$; its fiber over a
point $x\in M$ is the tangent space of $s^{-1}(x)\subset G$ at $1_{x}$. The target map $t: G \rightarrow M$
induces a so-called \textit{anchor map} $\rho: \gg \rightarrow TM$. Furthermore, the Lie bracket of vector fields
on $G$ induces a Lie bracket on the space $\Gamma\gg$ of sections. The anchor map respects this bracket and
induces a Lie algebra homomorphism $\rho: \Gamma \gg \rightarrow \mathfrak{X}(M)$ into the vector fields on $M$.
(There is no ``Lie's third theorem'' for Lie algebroids, stating that every such structure comes from some Lie
groupoid. For recent results on the integrability of Lie algebroids, see e.g. \cite{CraiFer, MoerMr}.)

\subsection{Some classes of Lie groupoids.}
Let $G$ be a Lie groupoid over $M$, with structure maps $s,t,m,u,i$ as above.
\begin{itemize}
\item[(a)] $G$ is called \textit{transitive} if the source and target together define a surjective submersion $(s,t): G \rightarrow M \times M$.
Transitive Lie groupoids are extensively discussed in \cite{Mack}.
\item[(b)] $G$ is called \textit{proper} if the pair $(s,t): G \rightarrow M \times M$ is a proper map (a closed map with compact Hausdorff
fibers).
\item[(c)] $G$ is called a \textit{bundle of Lie groups} if $s=t$. The bundle is called \textit{locally trivial} if $G$ is locally (in $M$)
of the form $K\times U$, where $K$ is a Lie group, $U\subset M$ is open and $s=t$ is the projection. A bundle of
compact Lie groups (i.e.  a bundle in which $s=t: G \rightarrow M$ is proper) is automatically locally trivial,
\cite{W}.
\item[(d)] $G$ is called a \textit{foliation groupoid} if the fibers of the pair $(s,t): G \rightarrow M\times M$ are discrete. These
groupoids are discussed further in \cite{CM, W}.
\item[(e)] An \textit{orbifold groupoid} is a proper foliation groupoid.
\item[(f)] $G$ is called an \textit{\'etale groupoid} if $s$ and $t$ are local diffeomorphisms (equivalently, if
$dim(G)=dim(M)$).
\item[(g)] $G$ is called a \textit{regular groupoid} if, for each $x\in M$, the target map $t: G\rightarrow M$ restricts to a map
$t:s^{-1}(x)\rightarrow M$ of constant rank (if $M$ is not connected, only locally of constant rank). This is
equivalent to requiring of the associated algebroid $\gg$ that its anchor map $\rho: \gg \rightarrow TM$ is a map
of vector bundles of locally constant rank. Any transitive Lie groupoid is regular, as is any \'etale groupoid and
any bundle of Lie groups.
\end{itemize}

\subsection{Foliations.}
A foliation (without singularities) on a manifold $M$ is an integrable subbundle $\FF$ of the tangent bundle $TM$.
One refers to the codimension of $\FF$ in $TM$ as the codimension of the foliation. Such a foliation gives rise to
a \textit{holonomy groupoid} $Hol(M,\FF)$, see \cite{Wi}. This is a foliation groupoid over $M$ whose associated
algebroid is the given bundle $\FF$, with the inclusion $\FF \hookrightarrow TM$ as the anchor map. There are in
general many other foliation groupoids over $M$ with this property, but $Hol(M,\FF)$ is in some sense the minimal
one (\cite{CM}, Proposition 1, page 183).

A \textit{leaf} of the foliation is a maximal connected immersed submanifold $L\subset M$ which is everywhere
tangent to $\FF$. (The codimension in $M$ of such a leaf is the codimension of $\FF$.) These leaves form a
partition of $M$.

A \textit{transversal} to a foliation $\FF$ of codimension $q$ is a $q$-dimensional submanifold $T\subset M$ which
is everywhere transversal to $\FF$. It is called a \textit{complete transversal} if $T$ intersects every leaf of
$\FF$. Complete transversals always exist (we do not require $T$ to be connected).

\subsection{Orbits.}
If $G$ is a regular groupoid over $M$, with algebroid $\gg$, the image of the anchor map $\rho: \gg \rightarrow
TM$ defines a foliation $\FF$ of $M$, called the \textit{orbit foliation} of $G$. If $x\in M$, the \textit{orbit}
of $G$ through $x$ is the subset $t s^{-1}(x)\subset M$. It is an immersed submanifold of $M$. The connected
components of the orbits of $G$ are precisely the leaves of the orbit foliation of $G$. One sometimes writes $M/G$
for the space of orbits of $G$, with the quotient topology.

\subsection{Homomorphisms and Morita equivalence.}
Let $G$ be a Lie groupoid over $M$ and $H$ be one over $N$. A homomorphism $\phi: G\rightarrow H$  consists of two
smooth maps (both denoted) $\phi: M \rightarrow N$ and $\phi: G\rightarrow H$, together commuting with all the
structure maps (i.e., $\phi(g_{1},g_{2})=\phi(g_{1})\phi(g_{2})$, $\phi(s g)=s(\phi g)$, etc). Such a $\phi$ is
called a (\textit{weak}) \textit{equivalence} if
\begin{itemize}
\item[(i)] $\phi$ induces a surjective submersion
\[t\pi_{1}:H\times_{N}M\rightarrow N\]
where $H\times_{N}M$ is the fibered product $\{(h,x)|h\in H, x\in X$, and $ s(h)=\phi(x)\}$;
\item[(ii)] The square
\[\xymatrix{G \ar[r]^{\phi} \ar[d]_{(s,t)}& H \ar[d]^{(s,t)} \\ M\times M \ar[r]^{\phi \times \phi} & N\times N}\]
is a fibered product.
\end{itemize}
Two groupoids $G_{1}$ over $M_{1}$ and $G_{2}$ over $M_{2}$ are said to be \textbf{Morita equivalent} if there is
a third groupoid $H$ (over $N$, say) for which there are weak equivalences $H\rightarrow G_{1}$ and $H \rightarrow
G_{2}$.

\subsection{Example.}
Let $G$ be a regular groupoid. Let $T$ be a complete transversal for the orbit foliation of $M$. Let $G_{T}=\{g\in
G | s(g),t(g)\in T\}$ be the restriction of $G$ to $T$. By transversality, $G_{T}$ is again a manifold, and
defines a Lie groupoid $G_{T}$ over $T$. The evident inclusion $G_{T} \hookrightarrow G$ is a weak equivalence, so
$G_{T}$ is Morita equivalent to $G$.

If $G$ is a foliation groupoid (1.3(d)) then $G_{T}$ is an \'etale groupoid(1.3(f)). So any foliation groupoid is
Morita equivalent to an \'etale groupoid (cf \cite{CM}).

If $G$ is a transitive groupoid, the one can take for $T$ a single point, so any transitive groupoid is Morita
equivalent to a Lie group.

If $G$ is a proper foliation groupoid, then $G_{T}$ is a proper \'etale groupoid. Such a groupoid determines a
unique orbifold structure on the space of orbits of $G$ (or of $G_{T}$), see \cite{MoerPronk}. In fact, one can
think of this orbifold as defined by the Morita equivalence class of $G$ or of $G_{T}$, whence the term ``orbifold
groupoid'' in 1.3(e) above.

\section{Regular groupoids as extensions of foliation groupoids}
In this section, $G$ is a fixed regular Lie groupoid over a manifold $M$.

\subsection{Isotropy.}
For a point $x\in M$, let $I(G)_{x}$ be the group of all arrows $x \rightarrow x$ in $G$. It is known \cite{Mack}
that $I(G)_{x}$ is submanifold of $G$ and a Lie group. These groups $I(G)_{x}$ for all $x$ together form a bundle
of Lie groups over M. This bundle will be denoted $I(G) \rightarrow M$, and referred to as the isotropy bundle of
$G$. For each $x\in M$, we write $I(G)_{x}^{\circ}$ for the connected component of the identity of the Lie group
$I(G)_{x}$. These groups $I(G)_{x}^{\circ}$ together form a bundle of Lie groups $I(G)^{\circ}$; we will exhibit
an explicit smooth structure in 2.5 below.

\subsection{Fiberwise components.}
We recall the following general facts. Let $f: N\rightarrow M$ be an arbitrary submersion between manifolds. Then
$f$ can be factored uniquely up to diffeomorphism as the composition
\[\xymatrix{N \ar[r]^{\pi} & C(f)\ar[r]^{e} & M}\]
where $C(f)$ is a (possibly non-Hausdorff) manifold, $e$ is a local diffeomorphism (i.e. an \'etale map), and
$\pi$ is a submersion with connected fibers. The fiber $e^{-1}(x)$ is discrete --- it is the set of connected
components of the manifold $f^{-1}(x)$. Any section $\alpha: M \rightarrow C(f)$ defines by pullback an open set
$N_{\alpha}=\{y\in N \,|\, \pi(y)=\alpha f(y)\}$ in $N$. In other words, the section $\alpha$ is a smoothly
varying choice of a component $\alpha(x)$ in $f^{-1}(x)$, and the union of these components defines an open
submanifold $N_{\alpha}$ of $N$.

\subsection{Restriction to a transversal.}\label{sec23}
Let $T\subset M$ be a complete transversal to the orbit foliation $\FF$ of $G$ on $M$, and let $G_{T}$ be the
restricted groupoid, as in 1.7, with source and target maps
\[\xymatrix{s,t: G_{T} \ar@<.5ex>[r] \ar@<-.5ex>[r] & M}.\]
Since $T$ is transversal to the orbits, this groupoid has the special property that the anchor map of its Lie
algebroid has rank zero. Or in other words, for each point $x\in T$ the target $t: G_{T} \rightarrow T$ restricts
to a locally constant map
\[t: s^{-1}(x) \rightarrow T.\]
A groupoid with this property is called \textit{essentially \'etale} in \cite{W}.

\subsection{The \'etale groupoid $E_{T}$.}
By 2.2 we can factor the submersion $s: G_{T} \rightarrow T$ as $s=\sigma\circ\pi$, where $\sigma:
E_{T}\rightarrow T$ is a local diffeomorphism and $\pi:G_{T}\rightarrow E_{T}$ is a submersion with connected
fibers. Since the target $t:G_{T}\rightarrow T$ is locally constant on the fibers of $s$, it factors uniquely as
$t=\tau\circ\pi$ for a map $\tau: E_{T}\rightarrow T$. In a similar fashion, the composition $m$ and inverse $i$
of $G_{T}$ induce unique composition and inverse maps on $E_{T}$. These make $E_{T}$ into an \'etale groupoid over
$T$, and $\pi: G_{T}\rightarrow E_{T}$ into a homomorphism of groupoids. The unit of the groupoid $E_{T}$ is the
composite $T \overset{u}{\to}G_{T}\overset{\pi}{\to} E_{T}$. This is a section of the \'etale map $\sigma:
E_{T}\rightarrow T$. As in 2.2, this section defines an open submanifold of $G_{T}$, whose fiber over a point
$x\in T$ is exactly the connected component $I(G_{T})_{x}^{\circ}=I(G)_{x}^{\circ}$ of the isotropy group at $x$.
Thus we obtain an ``extension'',
\begin{equation}\label{ext1}
\xymatrix{I(G_{T})^{\circ} \ar[r]^{j} & G_{T} \ar[r]^{\pi} & E_{T}},
\end{equation}
where $E_{T}$ is an \'etale groupoid over $T$, while $I(G_{T})^{\circ}$ is a bundle of connected Lie groups.
Furthermore, $j$ and $\pi$ are groupoid homomorphisms, $\pi$ is a submersion with connected fibers, and $j$ is an
open embedding.

This picture~(\ref{ext1}) extends to all of $M$, as in the following proposition.
\subsection{Proposition.}
\textit{Any regular Lie groupoid $G$ over $M$ fits into an extension
\begin{equation}\label{ext2}
\xymatrix{I(G)^{\circ}\ar[r]^{j} & G \ar[r]^{\pi} & E}
\end{equation}
where $I(G)^{\circ}$ is a bundle of connected Lie groups over $M$, $E$ is a foliation groupoid over $M$, $j$ and
$\pi$ are groupoid homomorphisms, $j$ is an embedding and $\pi$ an submersion with connected fibers.}

The properties of $j$ and $\pi$ determine $E$ and $I(G)^{\circ}$ uniquely up to isomorphism. The isotropy group
$I(E)_{x}$ of $E$ is the discrete group $I(G)_{x}/I(G)_{x}^{\circ}$ of connected components of $I(G)_{x}$. The
orbits of $E$ are the same as those of $G$, and define the same foliation $\FF$ on $M$. The restriction of $E$ to
$T$ is (isomorphic to) the \'etale groupoid $E_{T}$ in (\ref{ext1}) above.
\begin{proof}
It suffices to proof that $I(G)^{\circ}$ is a submanifold of $G$. For then we can view $I(G)^{\circ}$ as a
``normal'' (i.e. closed under conjugation by arrows) subgroupoid of $G$, and construct $E$ as the quotient
groupoid $E=G/I(G)^{\circ}$. Then $E$ is a manifold by the ``Godement Theorem'' for quotients, see \cite{Serre}.

By choosing local sections of appropriate submersions, we can find a cover of $M$ by open sets $U_{i}$ (foliation
charts for $\FF$) for which there are smooth maps $\gamma_{i}: U_{i}\rightarrow T$ and $\tau_{i}: U_{i}\rightarrow
G$ where $\tau_{i}(x): x\rightarrow \gamma_{i}(x)$. These maps give an isomorphism
\[G_{U_{i}}=U_{i}\times_{T}G_{T}\times_{T} U_{i},\]
mapping an arrow $g: x\rightarrow y$ in $G$ with $x,y\in U_{i}$ to the triple $(x, \tau_{i}(y)g \tau_{i}(x)^{-1},
y)$. Thus
\[I(G)^{\circ}|_{U_{i}}=U_{i}\times_{T}I(G_{T})^{\circ},\]
and this is a manifold because $I(G_{T})^{\circ}$ is (see 2.4) and $I(G_{T})^{\circ}\rightarrow T$ is a
submersion.

The proposition now follows.
\end{proof}
We single out the special case where the regular groupoid $G$ is proper. (Such groupoids are studied in [W].) In
this case, the groupoids $G_{T}$, $E_{T}$ and $E$ are all proper as well, and the bundles $I(G_{T})^{\circ}$ and
$I(G)^{\circ}$ are bundles of compact connected Lie groups, hence in particular locally trivial (cf 1.3). The
groupoid $E$ defines an orbifold structure on the orbit space $M/G=T/G_{T}$ (cf 1.7).
\subsection{Corollary.}
\textit{Any regular Lie groupoid $G$ over $M$ fits into an extension $K\rightarrow G\rightarrow E$ where $E$ is an
orbifold groupoid and $K$ is a locally trivial bundle of compact connected Lie groups.}

We will see later (e.g. 4.9 and 10.3) that such extensions are closely related to gerbes on orbifolds.

\subsection{General plan.}
Our strategy in this paper is to classify regular groupoids in two steps. First of all, given an extension of an
\'etale groupoid over $T$ of the form (\ref{ext1}), we describe all possible ways of expanding the extension to
one of the form (\ref{ext2}) over all of $M$. There is a cohomological obstruction to doing this (\S5). When the
obstruction vanishes, all these expansions to $M$ form a principal bundle under a group of bitorsors (\S4 and
\S6).

Secondly, we give a cohomological classification of all \'etale groupoids of the form (\ref{ext1}), (\S8). This
involves a new Cech cohomology of \'etale groupoids (\S7).

In this theory, non abelian cohomology plays a central r\^ole. For expository purposes, we treat the abelian case
first, and discuss non-abelian cohomology in \S9-12.

\newpage
\part{The principal action of bitorsors on regular groupoids}

\section{Bitorsors}
We review the basic definitions concerning bitorsors. Everything in this section essentially goes back to
\cite{Giraud}, although Giraud's presentation is in the context of sheaves and stacks on a site, while we are
concerned with bundles over a smooth manifold.

\subsection{K-torsors.}
Let $M$ be a manifold, and let $K$ be a bundle of Lie groups over $M$. A \textit{right} $K$\textit{-torsor} (or
\textit{principal} $K$\textit{-bundle}) is a manifold $P$, equipped with a surjective submersion $\pi:
P\rightarrow M$ and a smooth fiberwise action $P\times_{M} K \overset{\cdot}{\to} P$, with the property that the
induced map
\[(\pi_{1},\cdot):P\times_{M}K \rightarrow P\times_{M}P\]
is a diffeomorphism. Left $K$-torsors are defined in the obvious similar way.

Every right $K$-torsor $P$ determines a left $K$-torsor $P^{-1}$, with the same underlying space $P$. To
distinguish $P^{-1}$ from $P$, we denote the point $p\in P$ as $p^{-1}$ when thought of as a point in $P^{-1}$.
The left action of $K$ on $P^{-1}$ can then be written as
\[k\cdot p^{-1}=(p\cdot k^{-1})^{-1},\]
(where $x\in M$, $k\in K_{x}$, and $p\in P_{x}$ so $p^{-1}\in P_{x}^{-1}$).

Any right $K$-torsor $P$ gives rise to a \textit{division map}
\[\delta: P^{-1}\times_{M}P\rightarrow K\]
defined uniquely by the identity
\[q\cdot \delta(q^{-1},p)=p\]
In the sequel, we will simply write
\[\delta(q^{-1},p)=q^{-1} p.\]

\subsection{Bitorsors.}
Let $K$ and $L$ be two bundles of Lie groups over $M$. An $L$-$K$ \textit{bitorsor} is a manifold $P$ with a
surjective submersion $\pi:P\rightarrow M$, and structures of a left $L$-torsor as well as a right $K$-torsor,
which are compatible in the sense that
\[l (p k)= (l p) k\]
for any $l\in L$, $k\in K$ and $p\in P$, all over the same point in $M$. One sometimes writes
\[P=\, _{L}P_{K}\]
to emphasize that $P$ is an $L$-$K$ bitorsor.

\subsection{Tensor product.}
If $P$ is a right $K$-torsor and $Q$ is a left $K$-torsor, their tensor product $P\otimes_{K} Q$ is the manifold
obtained from the fibered product $P\times_{K}Q$ by identifying $(p k,q)$ and $(p, k q)$, for any $p\in P, q\in Q$
and $k\in K$ in the fibers over one and the same point in $M$. The equivalence class of $(p, q)$ is denoted by $p
\otimes q$. If $P=\, _{L}P_{K}$ and $Q=\, _{K}Q_{H}$ are bitorsors then $P\otimes Q$ carries obvious left-$L$ and
right-$H$ actions, making it a bitorsor. If $_{L}P_{K}$ is a bitorsor, then $P^{-1}$ is a $K$-$L$ bitorsor.

This tensor product is associative up to isomorphism. In fact, it can be viewed as composition in a category
$\mathsf{Bitor}(M)$, whose objects are bundles of Lie groups, and whose arrows $L\rightarrow K$ are isomorphism
classes of $L$-$K$ bitorsors. The identity arrow $K\rightarrow K$ in this category is represented by the trivial
torsor $K$. This category is in fact a groupoid, because if $P$ is an $L$-$K$ bitorsor, then the division map
$\delta: P^{-1}\times_{M}P \rightarrow K$ defines an isomorphism $\xymatrix{P^{-1}\otimes_{L} P \ar[r]^-{\sim} &
K}$; and similarly $P\otimes_{K}P^{-1}\cong L$.

\subsection{From torsors to bitorsors.}
Suppose $P$ is a right $K$-torsor over $M$, with associated division map $\delta: P^{-1}\times_{M}P\rightarrow K$,
$\delta(q^{-1},p)=q^{-1}p$. Then the manifold $P\otimes_{K} P^{-1}$ over $M$ carries the natural structure of a
bundle of Lie groups, with multiplication defined by the evident formula
\[(p\otimes q^{-1})(r \otimes s^{-1})= p(q^{-1} r)\otimes s^{-1}=p\otimes (q^{-1} r) s^{-1};\]
the unit section $M \rightarrow P\otimes_{K}P^{-1}$ is the map $x \mapsto p\otimes p^{-1}$ for any point $p\in
P_{x}$ (and independent of the choice of $p$). (This bundle is sometimes called the \textit{adjoint} or
\textit{conjugate} bundle of $P$. It is isomorphic to the isotropy bundle of the gauge groupoid of $P$; it can
also be constructed as the bundle obtained by twisting $K$ by the conjugate action of a cocycle defining $P$.)

There is a left action of this bundle of groups $P\otimes_{K} P^{-1}$ on $P$, defined in terms of the right action
of $K$ on $P$ by the obvious formula
\[(p\otimes q^{-1})\cdot r=p\cdot (q^{-1} r)\]
This action gives $P$ the structure of a left-$(P\otimes_{K}P^{-1})$ right-$K$ bitorsor.

More generally, for a bundle $L$ of Lie groups over $M$, $L$-$K$ bitorsor structures on $P$ correspond to
isomorphisms
\[\psi: P\otimes_{K}P^{-1}\rightarrow L.\]
Indeed, if $L$ acts on $P$ from the left, making $P$ into a bitorsor, then the division map $\varepsilon$ for this
left action, the map $\varepsilon: P\otimes_{K}P^{-1}\rightarrow L$ defined by the identity
$\varepsilon(p,q^{-1})\cdot q=p$, factors through the tensor product $P\otimes_{K}P^{-1}$, and gives an
isomorphism $\xymatrix{P\otimes_{K}P^{-1}\ar[r]^-{\sim}& L}$. Conversely, given such an isomorphism $\psi:
P\otimes_{K}P^{-1}\rightarrow L$, one can make $P$ into an $L$-$K$ bitorsor by defining the left action of $L$ on
$P$ by
\[l\cdot q=\mbox{the unique }q\mbox{ with }\psi(p,q^{-1})=l.\]
Thus, for any $L$-$K$ bitorsor $P$, the group $L$ and its action on $P$ are completely determined by $P$ and its
right $K$-torsor structure (cf. \cite{Giraud}, page 121).

\subsection{The abelian case.}
Let us assume that the bundle $K$ is a bundle of \textit{abelian} Lie groups. Let $P$ be a right $K$-torsor. Since
$K$ is abelian, the map
\[P\otimes_{K}P^{-1}\rightarrow K,\qquad p\otimes q^{-1}\rightarrow q^{-1}p\]
is well defined on the tensor product, and defines an isomorphism $P\otimes_{K}P^{-1}\simeq K$. Thus, by 3.4,
bitorsors correspond to automorphisms of the group bundle $K$. Explicitly, an isomorphism $\sigma: K\rightarrow K$
defines a left action of $K$ on $P$, by
\[k\cdot p=p\cdot \sigma(k).\]
Let us denote the bitorsor thus defined by $(P,\sigma)$.

For abelian $K$, there is also a familiar tensor product of two right $K$-torsors $P$ and $Q$. To distinguish it
from the one on bitorsors, let us temporarily denote it
\begin{equation}\label{ab:tens1}
P\otimes^{(r)}Q.
\end{equation}
It is the quotient of $P\times_{M}Q$ obtained by identifying $(p k,q)$ and $(p, q k)$.

The group $Aut(K)$ of automorphisms of the bundle $K$ acts from the \textit{left} on the right $K$-torsors: for
$\tau : K \overset{\sim}{\to} K$ and a right $K$-torsor $P$, there is a new right $K$-torsor $\tau_{!}(P)$ which
is the same space $P$ with new action $*$ defined by $p*k=p\cdot \tau^{-1}(k)$. For the tensor product of two
bitorsors $(P, \sigma)$ and $(Q,\tau)$ we then find the following formula
\begin{equation}\label{ab:tens2}
(P,\sigma)\otimes(Q,\tau)=(\tau_{!}(P)\otimes^{(r)}Q,\tau\sigma).
\end{equation}

Now let us pass to isomorphism classes. The tensor product~(\ref{ab:tens1}) makes the set of isomorphism classes
of right $K$-torsors into a group, denoted $H^{1}(M,K)$ as usual. Also, the tensor product 3.3 on bitorsors makes
the set of isomorphism classes of $K$-$K$ bitorsors into a group, which we denote by $Bitor(M,K)$. If we let
$Aut(K)^{op}$ be the opposite group, then this group acts on $H^{1}(M,K)$ from the right, and the
formula~(\ref{ab:tens2}) shows that for abelian $K$, there is an isomorphism of groups
\[Bitor(K)\simeq H^{1}(M,K)\rtimes Aut(K)^{op}.\]

\section{Action of bitorsors on regular groupoids}
We begin by taking up the discussion at the end of \S2. Let $G$ be a regular groupoid over a manifold $M$. Recall
that $G$ gives rise to an extension of the form
\begin{equation}\label{act:ext1}
K \overset{j}{\hookrightarrow} G \overset{\pi}{\twoheadrightarrow} E
\end{equation}
where $E$ is a foliation groupoid and $K$ is a bundle of (connected) Lie groups. For given $K$ and $E$, a map
between two such extensions $\xymatrix{K \ar[r]^{j} & G \ar[r]^{\pi} & E}$ and $\xymatrix{K \ar[r]^{j'} & G'
\ar[r]^{\pi'} & E}$ is a smooth homomorphism $\pi: G\rightarrow G'$ for which $\pi'\phi=\pi$ and $\phi j=j'$, as
usual. It follows that $\phi$ must be a diffeomorphism, hence an isomorphism of Lie groupoids. Thus, for $K$ and
$E$ fixed, there is a (discrete) groupoid of extensions of the form (\ref{act:ext1}), and maps between them. This
groupoid is denoted
\[Ext_{M}(E,K).\]

Let $T$ be a complete transversal to the orbit foliation of $E$ (or $G$) in (\ref{act:ext1}). Any extension of the
form (\ref{act:ext1}) can be restricted to $T$, to give an extension
\begin{equation}\label{act:ext2}
K_{T}\rightarrow G_{T}\rightarrow E_{T}
\end{equation}
over the manifold $T$. This operation defines a restriction functor $R$ between groupoids,
\begin{equation}\label{act:ext3}
R: Ext_{M}(E,K)\rightarrow Ext_{T}(E_{T},K_{T}).
\end{equation}
The purpose of this section is to show that, in some sense, this functor makes $Ext_{M}(E,K)$ into a principal
bundle over $Ext_{T}(E_{T},K_{T})$. The ``structure group'' of this bundle is formed by the $K$-$K$ bitorsors. A
precise formulation of the way this is a principal bundle involves some category theory, and will be given only in
\S6. In this section, however, we prove all the properties of extensions involved in the formulation, the main one
being Proposition 4.7.

\subsection{Twisting construction.}
Consider an extension of the form (\ref{act:ext1}), and an $L-K$ bitorsor $P$ over $M$, as in the previous
section. From these data one can form a new regular Lie groupoid over $M$, denoted $P \otimes_{K} G \otimes_{K}
P^{-1}$, which fits into an extension of the form
\[\xymatrix{L\ar[r]^-{j_{P}}& P \otimes_{K} G \otimes_{K} P^{-1} \ar[r]^-{\pi_{P}}& E},\]
as follows. For two points $x$ and $y$ in $M$, arrows $y \leftarrow x$ in $P \otimes_{K} G \otimes_{K} P^{-1}$ are
equivalence classes of triples
\[(q,g,p^{-1})\]
where $q\in P_{y}$, $g:x \rightarrow y$ in $G$, and $p\in P_{x}$. (Recall that $p^{-1}$ is just a notation for $p$
as an element of $P^{-1}_{x}$.) The equivalence relation is given by identifications $(q k, g, p^{-1})\sim (q, k
g, p^{-1})$ and $(q, g l, p^{-1})\sim (q, g, l p^{-1})$ (where $k\in K_{y}$, $l\in K_{x}$, and we have suppressed
$j: K\hookrightarrow G$ from the notation, thinking of $K$ as a subspace of $G$). The equivalence class of
$(q,g,p^{-1})$ is denoted
\[q\otimes g\otimes p^{-1}: x\rightarrow y.\]
For another such arrow
\[r\otimes h\otimes s^{-1}: y\rightarrow z,\]
the composition in the groupoid $P \otimes_{K} G \otimes_{K} P^{-1}$ is defined as
\[(r\otimes h\otimes s^{-1})(q\otimes g\otimes p^{-1}) =r\otimes h(s^{-1} q)g\otimes p^{-1}.\]
(Here $s^{-1} q\in K_{y}$ is the division $\delta(s^{-1},q)$ for the bitorsor $P$, defined by $s\cdot (s^{-1}
q)=q$, and we have simply written $h(s^{-1} q) g$ for the composition $h j(s^{-1} q) g$ in $G$, as before.) It is
easy to check that this composition is well defined on equivalence classes, and that the space of equivalence
classes is a manifold. Thus $P\otimes_{K}G\otimes_{K}P^{-1}$ is a Lie groupoid over $M$.

The map $\pi_{P}: P\otimes_{K}G\otimes_{K}P^{-1}\rightarrow E$,
\[\pi_{P}(q\otimes g\otimes p^{-1})=\pi(g),\]
is a groupoid homomorphism, and its kernel is precisely the adjoint bundle $P\otimes_{K}P^{-1}$ of Lie groups
considered in 3.4. This kernel is canonically isomorphic to $L$ (by the map $\varepsilon:
P\otimes_{K}P^{-1}\rightarrow L$ in 3.4). Thus we obtain an extension
\[\xymatrix{L\ar[r]^-{j_{P}}& P\otimes_{K}G\otimes_{K}P^{-1}\ar[r]^-{\pi_{P}}&E}\]
The map $j_{P}$ can explicitly be described for $l\in L_{x}$ by
\[j_{P}(l)=l q\otimes 1_{x}\otimes q^{-1},\]
independent of the choice of $q\in P_{x}$.

\subsection{Remark.}
Let us think of the twisting construction as a left action of $P$ on $G$ and denote it $P\cdot G$. Then it is
functorial, in the sense that if $_{H}Q_{L}$ is another bitorsor, there is a canonical isomorphism of extensions
\[Q\cdot (P\cdot G)= (Q\otimes_{L} P)\cdot G.\]

Consider an extension $K\rightarrow G\rightarrow E$ over $M$ and a complete transversal $T$ to the orbit foliation
$\FF$. We shall relate the twisting construction to the restriction functor (\ref{act:ext3}). For this we
introduce the following terminology.

\subsection{Sections of bitorsors.}
Let $P$ be an $L$-$K$ bitorsor over $M$. A section $\alpha: M\rightarrow P$ of $P$ induces an isomorphism of
bundles of Lie groups
\[\overline{\alpha}: K\rightarrow L\]
defined for $k\in K_{x}$ by $\alpha(x)k=\overline{\alpha}(k)\alpha(x)$. Or, using the notation $p q^{-1}$ for the
division map of $P$ for the left action by $L$,
\[\overline{\alpha}(k)=\alpha(x)k\alpha(x)^{-1}.\]
For this reason, we call $\overline{\alpha}$ the \textit{conjugation by} $\alpha$. If $L=K$ and $P$ is a $K-K$
bitorsor, we say that $\alpha$ is a \textit{central section} if $\overline{\alpha}$ is the identity. Note that an
arbitrary section $\alpha$ defines a trivialization $\xymatrix{K\ar[r]^-{\sim} &P}$ of right $K$-torsors by
$k\mapsto \alpha(x)k$, and that this is a trivialization of bitorsors iff $\alpha$ is central.

\subsection{Lemma}
\textit{Let $\alpha$ be a section over $T$ of the $L-K$ bitorsor $P$ over $M$. Then $\alpha$ induces an
isomorphism $\hat{\alpha}$ of extensions,
\[\xymatrix{K_{T}\ar[r]^{j}\ar[d]_{\overline{\alpha}}& G_{T}\ar[r]^{\pi}\ar[d]^{\hat{\alpha}}& E_{T}
\ar@{=}[d]\\L_{T}\ar[r]^-{j_{P}}& (P\otimes_{K}G\otimes_{K}P^{-1})_{T}\ar[r]^-{\pi_{P}}&E_{T}}.\] In particular,
if $L=K$ and $\alpha$ is central, then $\hat{\alpha}: G_{T}\rightarrow (P\otimes_{K}G\otimes_{K}P^{-1})_{T}$ is an
isomorphism in $Ext_{T}(E_{T}, K_{T})$.}

\begin{proof}
The map $\hat{\alpha}$ is defined, for an arrow $g: x\rightarrow y$ in $G_{T}$, by
\[\hat{\alpha}(g)=\alpha(y)\otimes g\otimes \alpha(x)^{-1}.\]
One checks that $\pi_{P}\circ\hat{\alpha}=\pi$ and $j_{P}\overline{\alpha}=\hat{\alpha}j$. (For the latter we have
for any $k\in K_{x}$,
\[\begin{array}{ll}
\hat{\alpha}j(k)&=\alpha(x)\otimes k\otimes \alpha(x)^{-1}\\
&=\alpha(x)k\otimes 1_{x}\otimes \alpha(x)^{-1}\\
&=\alpha(x)k\alpha(x)^{-1}\alpha(x)\otimes 1_{x}\otimes \alpha(x)^{-1}\\
&=\overline{\alpha}(k)\alpha(x)\otimes 1_{x}\otimes \alpha(x)^{-1}\\
&=j_{P}(\overline{\alpha}(k).)
\end{array}\]
\end{proof}

\subsection{Proposition.}
\textit{Let $P$ be an $L$-$K$ bitorsor over $M$ with a section $\alpha$ over $T$.
\begin{itemize}
\item[(i)]
Any extension of $\alpha$ to a section $\beta$ over $M$ gives an isomorphism $\hat{\beta}$ of extensions
\[\xymatrix{K\ar[r]^{j}\ar[d]_{\overline{\beta}}& G\ar[r]\ar[d]^{\hat{\beta}}& E
\ar@{=}[d]\\L\ar[r]^-{j_{P}}& P\otimes_{K}G\otimes_{K}P^{-1}\ar[r]&E}\] which restrict to $\hat{\alpha}$ over $T$.
\item[(ii)]
If $\xymatrix{v:K\ar[r]^-{\sim}&L}$ and $\xymatrix{u:G\ar[r]^-{\sim}&P\otimes_{K}G\otimes_{K}P^{-1}}$ form an
isomorphism of extensions
\[\xymatrix{K\ar[r]^{j}\ar[d]_{v}& G\ar[r]\ar[d]^{u}& E
\ar@{=}[d]\\L\ar[r]^-{j_{P}}& P\otimes_{K}G\otimes_{K}P^{-1}\ar[r]&E}\] which restricts to $\hat{\alpha}$ over
$T$, then $u=\hat{\beta}$ and $v=\overline{\beta}$ for a unique section $\beta$ over $M$, which restricts to
$\alpha$ over $T$.
\item[(iii)]
In particular, if $K=L$ and $v$ is the identity, then $\beta$ is central whenever $\alpha$ is.
\end{itemize}}

\begin{proof}
Part (i) is proved just as in 4.4, and (iii) follows immediately from (ii). We prove (ii). First of all, by
choosing local sections of the surjective submersion $s: G\times_{M}T=t^{-1}(T)\rightarrow M$, we find an open
cover $\{U_{i}\}$ of $M$, together with mappings $\gamma_{i}: U_{i}\rightarrow T$ and $\tau_{i}:U_{i}\rightarrow
G$, such that $\gamma_{i}$ is a submersion and $\tau_{i}(x):x\rightarrow \gamma_{i}(x)$ in $G$. Write $\tau_{i
j}:U_{i j}\rightarrow G$ for the map $\tau_{i j}(x)=\tau_{i}(x)\tau_{j}^{-1}(x):\gamma_{j}(x)\rightarrow
\gamma_{i}(x)$ in $G_{T}$, where $U_{ij}=U_{i}\cap U_{j}$ as usual. To define $\beta$, choose for each $x\in M$ an
$i$ with $x\in U_{i}$. Then $u(\tau_{i}(x))$ can be written in the form
\[u(\tau_{i}(x))=\alpha(\gamma_{i}(x))\otimes\tau_{i}(x)\otimes p^{-1}\]
for a unique $p=p_{i}(x)\in P_{x}$. This $p$ does not depend on the choice of $i$. Indeed, if $x\in U_{i j}$, then
$\tau_{i j}(x)\tau_{j}(x)=\tau_{i}(x)$ in $G$, so $u(\tau_{i j}(x))u(\tau_{j}(x))=u(\tau_{i}(x))$ in
$P\otimes_{K}G\otimes_{K}P^{-1}$. Since $u$ extends $\hat{\alpha}$ by assumption, this means that
\[\begin{array}{c}(\alpha(\gamma_{i}(x))\otimes\tau_{i j}(x)\otimes\alpha(\gamma_{j}(x))^{-1})\circ
(\alpha(\gamma_{j}(x))\otimes\tau_{j}(x)\otimes p_{j}(x)^{-1})\\
=\alpha(\gamma_{i}(x))\otimes\tau_{i}(x) \otimes p_{i}(x)^{-1},
\end{array}\] so
\[\alpha(\gamma_{i}(x))\otimes\tau_{i j}(x)\tau_{j}(x)\otimes p_{j}(x)^{-1}=\alpha(\gamma_{i}(x))\otimes
\tau_{i}(x)\otimes p_{i}(x)^{-1}.\] Since $\tau_{i j}(x)\tau_{j}(x)=\tau_{i}(x)$, we conclude $p_{i}(x)=p_{j}(x)$.
Thus we can define
\[\beta(x)=p_{i}(x) \mbox{  for }x\in U_{i},\]
independent of $i$. Then $\hat{\beta}=u$ and $\beta$ extends $\alpha$. The section $\beta$ is unique, because
$\beta(x)=p_{i}(x)$ is completely determined by the equation
\[u(\tau_{i}(x))=\beta(\gamma_{i}(x))\otimes\tau_{i}(x)\otimes\beta(x)^{-1}\]
together with the fact that $\beta(\gamma_{i}(x))=\alpha(x)$. Finally, since $\hat{\beta}=u$, we have $\hat{\beta}
j= j_{P} v$, while by (i) $\hat{\beta}j= j_{P}\overline{\beta}$. Thus $v=\overline{\beta}$.
\end{proof}

\subsection{Remark.}
For later use, we observe that the groupoid $G$ over $M$ can be reconstructed from its restriction $G_{T}$ and the
``cocycle'' given by the maps $\gamma_{i}: U_{i}\rightarrow T$ and $\tau_{i j}(x):\gamma_{j}(x)\rightarrow
\gamma_{i}(x)$, as in the previous proof. Indeed, arrows $g:x\rightarrow y$ in $G$ correspond to equivalence
classes of triples $(i,h,j)$ where $x\in U_{i}$, $y\in U_{j}$ and $h:\gamma_{i}(x)\rightarrow \gamma_{j}(x)$ (via
the correspondence $h\tau_{i}(x)=\tau_{j}(y)g$). Two such triples $(i,h,j)$ and $(i',h',j')$ are equivalent if
\[\tau_{j' j}(y) h=h'\tau_{i' i}(x).\]
If we are just given the groupoid $G_{T}$ and cocycle $\tau=\{\gamma_{i},\tau_{i j}\}$, we can construct an
abstract groupoid $G=G_{T}$ over $M$ in this way. It is a Lie groupoid because the $\gamma_{i}$ are submersions,
and its restriction to $T$ is isomorphic to $G_{T}$.

\subsection{Proposition.}
\textit{Let $K\overset{i}{\hookrightarrow}G\overset{\pi}{\rightarrow}E$ and $K\overset{j}{\hookrightarrow}H
\overset{\rho}{\rightarrow}E$ be two extensions over $M$, and suppose $\phi: H_{T}\rightarrow G_{T}$ is an
isomorphism between the restricted extensions $K_{T}\rightarrow G_{T}\rightarrow E_{T}$ and $K_{T}\rightarrow
H_{T}\rightarrow E_{T}$. Then there exists a $K-K$ bitorsor $P$ over $M$ with a central section $\alpha$ over $T$,
for which there is an isomorphism
\[\psi: P\otimes_{K}H\otimes_{K}P^{-1}\rightarrow G\]
with the property that
\[\psi_{T}\circ\hat{\alpha}=\phi:H_{T}\rightarrow G_{T}\]}

\begin{proof}
Consider the manifold $Q$ of pairs of arrows
\[x \overset{g}{\gets}t\overset{h}{\gets}x\]
where $x\in M$, $t\in T$ is a point in the transversal, and $g$ and $h$ are arrows in $G$ and $H$ respectively
such that $\pi(g)\rho(h)=1_{x}$ in $E$. ($Q$ is indeed a manifold because it can be constructed as the pullback of
the unit $M\rightarrow E$ along the submersion $\xymatrix{G\times_{M}T\times_{M}H\ar[r]^-{\pi\times
1\times\rho}&E\times_{M}T\times_{M}E\ar[r]^-{m}&E}$.) The groupoid $H_{T}$ acts from the right on $Q$ along the
map $Q\rightarrow T$ sending $(x \overset{g}{\gets}t\overset{h}{\gets}y)$ to $t$, by
\[(x\overset{g}{\gets}t\overset{h}{\gets}x)\cdot(t\overset{f}{\gets}t')=(x\overset{g\phi(f)}{\gets}t'\overset{f^{-1} h}{\gets}x)\]
This is a principal action, and we define $P$ to be the quotient. Thus points in $P$ can be written
\[g\otimes h\]
where $g: t\rightarrow x$ in $G$, $h: x\rightarrow t$ in $H$; and $g \phi(f)\otimes h=g\otimes f h$ for any
$f:t'\rightarrow t$. This manifold $P$ is a bitorsor under the actions by $K$,
\[\begin{array}{rl}
k\cdot(g\otimes h)&=i(k)g\otimes h\\
(g\otimes h)\cdot k&=g\otimes h j(k),
\end{array}\]
for $k\in K_{x}$. Indeed, both actions are principal. (For example, to see that the left one is, consider two
pairs $x\overset{g_{1}}{\gets}t_{1}\overset{h_{1}}{\gets}x$ and
$x\overset{g_{2}}{\gets}t_{2}\overset{h_{2}}{\gets}x$ representing points in $P_{x}$. Write
$f=h_{1}h_{2}^{-1}:t_{2}\rightarrow t_{1}$. Then $g_{1}\otimes h_{1}=g_{1}\otimes f h_{2}=g_{1}\phi(f)\otimes
h_{2}$, and $k=g_{2}(g_{1} \phi(f))^{-1}$ is the unique one with $k\cdot(g_{1}\otimes h_{1})=g_{2}\otimes h_{2}$.)
This bitorsor $P$ has a canonical section $\alpha$ over $T$, defined by
\[\alpha(t)=1_{t}\otimes 1_{t}\]
(the left $1_{t}$ being the identity at $t$ in $G$, the right $1_{t}$ the one in $H$). This section is obviously
central. The isomorphism
\[\psi:P\otimes_{K}H\otimes_{K}P^{-1}\rightarrow G\]
is the map defined, for $h:x\rightarrow y$ in $H$ and $x\overset{g_{1}}{\gets}t_{1}\overset{h_{1}}{\gets}x$
representing an element in $P_{x}$ and $y\overset{g_{2}}{\gets}t_{2}\overset{h_{2}}{\gets}y$ one in $P_{y}$, by
\[\psi((g_{2}\otimes h_{2})\otimes h\otimes (g_{1}\otimes h_{1}))=g_{2}\circ \phi(h_{2} h h_{1}^{-1})\circ g_{1}^{-1}.\]
One checks that $\pi\circ\psi=\rho_{P}:P\otimes_{K}H\otimes_{K}P^{-1}\rightarrow E$, and that $\psi\circ j_{P}=i:
K\rightarrow  G$. Note, finally, that for $h: x\rightarrow y$ in $H_{T}$, we have
\[\begin{array}{ll}
\psi(\hat{\alpha}(h))&=\psi((1\otimes 1)\otimes h\otimes (1\otimes 1)^{-1})\\
&=\phi(h)
\end{array}\]
Thus $\psi \hat{\alpha}=\phi$, as asserted in the proposition.
\end{proof}

\subsection{Remark.}
We consider briefly the case where $K$ is abelian. If $K\overset{j}{\to}G\overset{\pi}{\to}E$ is an extension and
$K$ is abelian, then $E$ acts on $K$ in the usual way: for an arrow $e: x\rightarrow y$ in $E$ and any $k\in
K_{y}$, the right action $k\cdot e$ is defined by
\[j(k\cdot e)=g^{-1} j(k) g\]
where $g:x\rightarrow y$ is any arrow in $G$ with $\pi(g)=e$. Then $k\cdot e$ is independent of the choice of $g$
because $K$ is abelian, and the action is smooth because $\pi:G\rightarrow E$ is a submersion.

By 3.5, the $K-K$ bitorsors are of the form $(P, \sigma)$, where $P$ is a right $K$-torsor and $\sigma\in Aut(K)$
defines the left action of $K$ on $P$, via $k\cdot p=p\cdot \sigma(k)$. For such a $K-K$ torsor, the twisted
extension $K\rightarrow (P,\sigma)\otimes_{K}G\otimes_{K}(P,\sigma)^{-1} \rightarrow E$ induces a new action of
$E$ on $K$. A small calculation shows that this new action is obtained by conjugation of the old one by $\sigma$.
In particular, the new action coincides with the old one iff $\sigma$ is $E$-equivariant; i.e. $\sigma_{y}(k)\cdot
e= \sigma_{x}(k\cdot e)$ for any $e:x\rightarrow y$ in $E$ and $k\in K_{y}$.

Furthermore, if $(P,\sigma)$ possesses a central section $\alpha$ over $T$, then $\sigma=id$ on $T$. If $\sigma$
is moreover $E$-equivariant, then this implies that $\sigma=id$ on all of $M$.

Thus, we see that if $K\rightarrow G\rightarrow E$ and $K\rightarrow H\rightarrow E$ are two extensions inducing
the same action of $E$ on $K$, and if $\phi: H_{T}\rightarrow G_{T}$ is an isomorphism between the restricted
extensions over $T$, then the bitorsor of Proposition 4.7 is of the form $(P, id)$. Any section of such a bitorsor
is automatically central.

\subsection{Example.}
Suppose the bundle $K$ of abelian Lie groups over $M$ is trivial, i.e. of the form $M\times A$ where $A$ is a
fixed abelian Lie group. An extension
\[M\times A\overset{j}{\to}G\overset{\pi}{\to}E\]
of Lie groupoids over $M$ is called a \textit{central} extension of $E$ by $A$ if the induced action of $E$ on
$M\times A$ is trivial. If we write $j(x,a)=a\in G_{x}$, this simply means that $a g=g a$ for any $g:x\to y$ in
$G$ and any $a\in A$. Such central extensions are closely related to bundle gerbes. In fact, a bundle gerbe in the
sense of Murray \cite{Mu} is the same thing as a central extension by $S^{1}$ of a pair groupoid
$\xymatrix{X\times_{M}X \ar@<.5ex>[r] \ar@<-.5ex>[r]& M}$ given by a surjective submersion $X\to M$. Bundle gerbes
over orbifolds \cite{LupUrb} are central extensions by $S^{1}$ of proper \'etale groupoids $E$.

\section{The obstruction class}
We now consider the question whether a given extension $K_{T}\rightarrow G_{T}\rightarrow E_{T}$ is isomorphic to
the restriction of an extension $K\rightarrow G\rightarrow E$; in other words, whether the functor (9) in \S4 is
surjective up to isomorphism. For the moment, we assume that the bundle $K$ is a bundle \textit{abelian} Lie
groups. The general case involves non-abelian cohomology and will be postponed until \S10.

If $\xymatrix{K\ar[r]^{j}&G\ar[r]^{\pi}&E} $ is an extension and $K$ is abelian, then $E$ acts on $K$ as usual
(see 4.8). Thus, $K$ is an $E$-equivariant bundle of abelian Lie groups. Similarly $K_{T}$ is an
$E_{T}$-equivariant bundle. In fact, since $E_{T}$ and $E$ are Morita equivalent Lie groupoids, the restriction
operation defines an equivalence of categories between $E$-equivariant bundles of abelian Lie groups over $M$ and
$E_{T}$-equivariant such bundles over $T$.

Now fix a foliation groupoid (or, for that matter, any regular Lie groupoid) $E$ over $M$, and let $E_{T}$ be its
restriction to the transversal $T$. Let $K_{T}\rightarrow B\rightarrow E_{T}$ be an extension of $E_{T}$ by a
bundle of abelian Lie groups, as before. Thus $B$ is a Lie groupoid over $T$, and $K_{T}$ is a $E_{T}$-equivariant
bundle over $T$. Let $K$ be the corresponding $E$-equivariant bundle over $M$. Let $\underline{K}$ be the sheaf on
$M$ of smooth local sections of $K$.

\subsection{Proposition.} \textit{There is an obstruction $\lambda(B)\in
\check{H}^{2}(M,\underline{K})$ to expanding $K_{T}\rightarrow B \rightarrow E_{T}$ to an extension $K\rightarrow
G\rightarrow E$ over $M$.}
\subsection{Corollary.} \textit{If $\check{H}^{2}(M,\underline{K})=0$ then the
restriction functor (cf. \S4(9))
\[R:Ext_{M}(E,K)\rightarrow Ext_{T}(K_{T},E_{T})\]
is surjective up to isomorphism.}

The Proposition and its Corollary also hold in the non-abelian case, and are stated and proved in this generality
in \S10.

In the proof of Proposition 5.1, we shall make the additional assumption that the bundle $K$ is locally trivial.
This assumption is always fulfilled if $K$ is a bundle of compact Lie groups, cf. \S1. More importantly, the
assumption isn't really necessary, but in the more general case the obstruction class does not lie in the \v{C}ech
cohomology $\check{H}^{2}(M,\underline{K})$ but in the cohomology for a hypercover; cf. also the discussion in
\cite{Breen}, page 38.

\begin{proof} (of Proposition 5.1)
Consider an extension
\[\xymatrix{K_{T}\ar[r]& B\ar[r]^{\pi}& E_{T}},\]
where $E_{T}$ and $K_{T}$ are the restrictions of $E$ and $K$, as introduced just before the statement of the
proposition. As before, we will delete the embedding $K_{T}\hookrightarrow B$ from the notation, and view $K_{T}$
as a subset of $B$. As in the proof of Proposition 4.4(ii), there exists an open cover $\{U_{i}\}$ of $M$ and
submersions $\gamma_{i}: U_{i}\rightarrow T$, together with mappings $\tau_{i}: U_{i}\rightarrow E$ such that
$\tau_{i}(x):x\rightarrow \gamma_{i}(x)$ in $E$. These $\tau_{i}$ define a cocycle $\tau_{i j}: U_{i j}\rightarrow
E_{T}$
\[\tau_{i j}(x)=\tau_{i}(x)\tau_{j}(x)^{-1}:\gamma_{j}(x)\rightarrow \gamma_{i}(x)\]
and $E$ can be reconstructed from $E_{T}$ and this cocycle (cf. Remark 4.6). If we can lift this cocycle to a
cocycle $\hat{\beta}=\{\hat{\beta}_{i j}\}$,
\[\hat{\beta}_{i j}(x):\gamma_{j}(x)\rightarrow \gamma_{i}(x) \mbox{ in }B\]
which is mapped by $\pi: B\rightarrow E_{T}$ to the cocycle $\{\tau_{i j}\}$, then this cocycle $\hat{\beta}$ can
be used to define a Lie groupoid $G=G_{\tau}$ over $M$, which fits into an extension $K\rightarrow G\rightarrow E$
restricting to $K_{T}\rightarrow B\rightarrow E_{T}$, up to isomorphism. The cohomology class $\lambda(B)$ will be
constructed as an obstruction to such a lifting $\{\hat{\beta}_{i j}\}$ of $\{\tau_{i j}\}$.

First of all, since $\xymatrix{B\ar[r]^{\pi}&E_{T}}$ is a surjective submersion, we can find for each $i$ and $j$
an open cover
\[U_{i j}=\bigcup V_{\xi}\]
of $U_{i j}$ by open sets $V_{\xi}$, over which there exists a lifting $\beta_{\xi}: V_{\xi}\rightarrow B$ of
$\tau_{i j}$,
\[\beta_{\xi}(x):\gamma_{j}(x)\rightarrow \gamma_{i}(x) \mbox{ in }B \qquad (x\in V_{\xi}).\]
Then for two indices $\xi$ and $\zeta$ and $x\in V_{\xi \zeta}$, the arrow
$\beta_{\xi}(x)\beta_{\zeta}(x)^{-1}:\gamma_{i}(x)\rightarrow \gamma_{i}(x)$ lies in the kernel of $\pi:
B\rightarrow E_{T}$, hence defines a cocycle $c=\{c_{\xi \zeta}\}$ on $U_{i j}$ with values in the sheaf
$\underline{\gamma_{i}^{*}(K_{T})}$ of sections of $\gamma_{i}^{*}(K_{T})$, by
\[c_{\xi \zeta}(x)=\beta_{\xi}(x)\beta_{\zeta}(x)^{-1}.\]
Now, if we choose the original cover $\{V_{\xi}\}$ of $U_{i j}$ to be ``good'' or ``Leray'' covers, then
$\check{H}^{1}(U_{i j},\underline{\gamma_{i}^{*}(K_{T})})=0$. Thus $c=\{c_{\xi, \zeta}\}$ is a coboundary, say
$c=d f$ where $f_{\xi}\in\Gamma(V_{\xi},\underline{\gamma_{i}^{*}(K_{T})})$. Then $f_{\xi}^{-1}f_{\zeta}=c_{\xi
\zeta}=\beta_{\xi}(x)\beta_{\zeta}(x)^{-1}$ for each $x\in V_{\xi}$, so the maps $f_{\xi}(x)\cdot \beta_{\xi}(x):
V_{\xi}\rightarrow B$ agree on overlaps, and define a map $\beta_{i j}: U_{i j}\rightarrow B$,
\[\beta_{i j}(x)=f_{\xi}(x)\beta_{\xi}(x):\gamma_{j}(x)\rightarrow \gamma_{i}(x)\]
(where $x\in V_{\xi}$, independent of $\xi$), which maps to $\tau_{i j}$.

Now define a 2-cocycle $\nu=\nu_{i j k}$ by setting, for $x\in U_{i j k}=U_{i}\cap U_{j}\cap U_{k}$,
\[\nu_{i j k}(x)=(\beta_{i j}(x)\beta_{j k}(x)\beta_{i k}(x)^{-1})\cdot \tau_{i}(x)\]
Here the composition $\beta_{i j}(x)\beta_{j k}(x)\beta_{i k}(x)^{-1}$ lies in the kernel of $\pi$, i.e. in
$K_{\gamma_{i}(x)}$, and the $\tau_{i}(x)$ at the right refers to the action of $E$ on $K$, so that
\[\nu_{i j k}(x)\in K_{x}.\]
Using that the action of $\tau_{i}(x)$ is defined by conjugation by any arrow $g:x\rightarrow \gamma_{i}(x)$ with
$\pi(g)=\tau_{i}(x)$, independent of the choice of $g$, it is easily verified that $\{\nu_{i j k}\}$ satisfies the
cocycle condition $\nu_{j k l}-\nu_{i k l}+\nu_{i j l}-\nu_{i j k}$=0. Thus $\{\nu_{i j k}\}$ defines a \v{C}ech
cohomology class
\[[\nu]\in \check{H}^{2}(M, \underline{K}).\]
It is now straight forward to verify the following two properties, as we will see shortly:
\begin{itemize}
\item[(i)]
The class $[\nu]$ does not depend on the choice of the $\{\beta_{i j}\}$ lifting $\{\tau_{i j}\}$.
\item[(ii)]
If $[\nu]=0$ then the $\beta_{i j}$ can be modified into a cocycle $\beta_{i j}'$ lifting $\tau_{i j}$.
\end{itemize}
As said, such a cocycle $\{\beta_{i j}'\}$ allows us to construct an extension $K\rightarrow G\rightarrow E$ of
groupoids over $M$, the restriction of which to $T$ is isomorphic to $K_{T}\rightarrow B\rightarrow E_{T}$. Thus,
if we define $\lambda(B)=[\nu]$, the proposition is proved. To conclude, then, we verify (i) and (ii).

For (i), suppose $\{\alpha_{i j}\}$ is another lift of $\{\tau_{i j}\}$. Then $\alpha_{i j}$ and $\beta_{i
j}:\gamma_{j}\rightarrow \gamma_{i}$ have the same image $\tau_{i j}(x)$ in $E_{T}$, so we can write $\beta_{i
j}(x)=\psi_{i j}(x)\alpha_{i j}(x)$ for a unique arrow $\psi_{i j}(x):\gamma_{i}(x)\rightarrow \gamma_{i}(x)$ in
$K_{\gamma_{i}(x)}$. Write
\[\tilde{\psi}_{i j}(x)=\psi_{i j}(x)\cdot \tau_{i}(x)\in K_{x}\qquad(x\in U_{i j}).\]
This defines a cochain $\{\tilde{\psi}_{i j}\}$ with values in $\underline{K}$. For this cochain we have (deleting
the argument $x\in U_{i j k}$ from the notation)
\[\begin{array}{ll}
(d \tilde{\psi})_{i j k}&=(\psi_{j k}\cdot \tau_{j})(\psi_{i k}^{-1}\cdot\tau_{i})(\psi_{i j}\cdot\tau_{i})\\
&=((\psi_{j k}\cdot \tau_{i j}^{-1})\psi_{i k}^{-1}\psi_{i j})\cdot\tau_{i}.
\end{array}\]
Thus
\[(d\tilde{\psi})_{i j k}((\alpha_{i j}\alpha_{j k}\alpha_{i k}^{-1})\cdot\tau_{i})
=((\psi_{j k}\cdot \tau_{i j}^{-1})\psi_{i k}^{-1}\psi_{i j}\alpha_{i j}\alpha_{j k}\alpha_{i k}^{-1})\cdot
\tau_{i}\] and one calculates,
\[\begin{array}{lr}
(\psi_{j k}\cdot \tau_{i j}^{-1})\psi_{i k}^{-1}\psi_{i j}\alpha_{i j}\alpha_{j k}\alpha_{i k}^{-1}=\\
\qquad=(\psi_{j k}\cdot \tau_{i j}^{-1})\psi_{i k}^{-1}\beta_{i j}\psi_{j k}^{-1}\beta_{j k}\beta_{i
k}^{-1}\psi_{i k}
&\mbox{(by def. of $\psi$)}\\
\qquad=(\psi_{j k}\cdot \tau_{i j}^{-1})\beta_{i j}\psi_{j k}^{-1}\beta_{j k}\beta_{i k}^{-1}
&\mbox{(because }K\mbox{ is abelian)}\\
\qquad=(\beta_{i j}\psi_{j k})\psi_{j k}^{-1}\beta_{j k}\beta_{i k}^{-1}
&\mbox{(action by }\tau_{i j}\mbox{ is conjugation by }\beta_{i j})\\
\qquad=\beta_{i j}\beta_{j k}\beta_{i k}^{-1}
\end{array}\]
This shows that the 2-cocycles $\{(\beta_{i j}\beta_{j k}\beta_{i k}^{-1})\cdot\tau_{i}\}$ and $\{(\alpha_{i
j}\alpha_{j k}\alpha_{i k}^{-1})\cdot\tau_{i}\}$ indeed differ by a coboundary $d\tilde{\psi}$ with value in
$\underline{K}$.

To conclude the proof, we verify (ii). Suppose $\nu=d \mu$ for a 2-cochain $\mu=\{\mu_{i j}\}$ with values in
$\underline{K}$. Thus
\begin{equation}\label{obs:eq1}
(\beta_{i j}\beta_{j k}\beta_{i k}^{-1})\cdot\tau_{i}=\mu_{j k}\mu_{i k}^{-1}\mu_{i j}.
\end{equation}
Then we can replace the lift $\beta_{i j}$ of $\tau_{i j}$ by $\hat{\beta_{i j}}$,
\[\hat{\beta_{i j}}(x)=\beta_{i j}(x)(\mu_{i j}(x)^{-1}\cdot\tau_{j}(x)^{-1}).\]
These new $\hat{\beta_{i j}}$ indeed form a cocycle: first rewrite (10) as
\begin{equation}\label{obs:eq2}
\begin{array}{ll}
\beta_{i j}\beta_{j k}&=((\mu_{j k}\mu_{i k}^{-1}\mu_{i j})\cdot\tau_{i}^{-1})\beta_{i k}\\
&=\beta_{i k}\beta_{i k}^{-1}((\mu_{j k}\mu_{i k}^{-1}\mu_{i j})\tau_{i}^{-1})\beta_{i k}\\
&=\beta_{i k}((\mu_{j k}\mu_{i k}^{-1}\mu_{i j})\cdot\tau_{k}^{-1}).
\end{array}
\end{equation}
Then
\[\begin{array}{llr}
\hat{\beta_{i j}}\hat{\beta_{j k}}&=\beta_{i j}(\mu_{i j}^{-1}\cdot\tau_{j}^{-1})\beta_{j k}(\mu_{j k}^{-1}\cdot \tau_{k}^{-1})\\
&=\beta_{i j}\beta_{j k}\beta_{j k}^{-1}(\mu_{i j}^{-1}\cdot\tau_{j}^{-1})\beta_{j k}(\mu_{j k}^{-1}\cdot \tau_{k}^{-1})\\
&=\beta_{i j}\beta_{j k}(\mu_{i j}^{-1}\cdot\tau_{j}^{-1}\cdot\tau_{j k})(\mu_{j k}^{-1}\cdot\tau_{k}^{-1})\\
&=\beta_{i j}\beta_{j k}(\mu_{i j}^{-1}\mu_{j k}^{-1})\cdot\tau_{k}^{-1})\\
&=\beta_{i k}(\mu_{i k}\cdot\tau_{k}^{-1})&\mbox{(by (\ref{obs:eq2}))}\\
&=\hat{\beta}_{i k}.
\end{array}\]
\end{proof}

\section{Monoidal categories and principal bundles}

The results of the previous sections express that the restriction functor from extensions $K\rightarrow
G\rightarrow E$ over $M$ to such extensions over $T$ is a principal bundle for the ``group'', given by the tensor
product on $K-K$ bitorsors over $M$ equipped with a central section over $T$. In this section, we provide the
categorical framework for a precise formulation of this result.

\subsection{Actions of monoidal categories.}
Recall from \cite{CWM} that a \textit{monoidal} (or \textit{tensor}) category is a category $\PP$ equipped with a
functorial tensor product $\otimes$, which is associative up to coherent natural isomorphism $a=a_{X, Y, Z}$,
\[a:(X\otimes Y)\otimes Z\overset{\sim}{\rightarrow}X\otimes(Y\otimes Z)\]
(here $X$, $Y$ and $Z$ are objects of $\PP$), and has a two-sided unit $I$, as expressed by coherent natural
isomorphisms $r=r_{X}$ and $l=l_{X}$
\[r:X\otimes I\overset{\sim}{\rightarrow}X \mbox{ and }l: I\otimes X\overset{\sim}{\rightarrow}X,\]
for each object in $\PP$. An \textit{action} of such a monoidal category $\PP$ on a category $\EE$ is given by a
functor
\[\PP\times\EE\rightarrow \EE\]
denoted $(P, E)\mapsto P\cdot E$, together with natural isomorphisms $u=u_{E}$ and $t=t_{X,Y,E}$,
\[\begin{array}{l}
u:I\cdot E \overset{\sim}{\rightarrow} E\\
t:X\cdot(Y\cdot E)\overset{\sim}{\rightarrow}(X\otimes Y)\cdot E,
\end{array}\]
which are compatible with the structural isomorphisms $a$, $l$ and $r$ of $\PP$, in the obvious sense. For
example, for $X,Y,Z\in \PP$ and $E\in\EE$, the pentagon
\[\xymatrix{X\cdot(Y\cdot(Z\cdot E))\ar[d]_{t}\ar[rr]^-{t}&&(X\otimes Y)\cdot(Z\cdot E)\ar[d]^{t}\\
X\cdot((Y\otimes Z)\cdot E)\ar[dr]_{t}&&((X\otimes Y)\otimes Z)\cdot E\ar[dl]_{a}\\
&(X\otimes(Y\otimes Z))\cdot E}\] should commute.

If $R:\EE\rightarrow \BB$ is a functor, a \textit{fiberwise action} of $\PP$ on $\EE$ is given by an action of
$\PP$ on $\EE$ as above, together with natural isomorphisms $b=b_{X, E}$ in $\BB$,
\[b: R(X\cdot E)\overset{\sim}{\rightarrow} R(E),\]
compatible with the earlier structure. For example, the square
\[\xymatrix{R(X\cdot(Y\cdot E))\ar[r]^{R t}\ar[d]_{b}&R((X\otimes Y)\cdot E)\ar[d]^{b}\\
R(Y\cdot E)\ar[r]^{b}&R(E)}\] should commute, and $b_{I, E}=b: R(I\cdot E)\rightarrow R(E)$ should coincide with
$R(u)$.

The (weak) categorical pullback $\EE\times_{\BB}\EE$ is the category of triples $(E,F,\alpha)$, where $E$ and $F$
are objects of $\EE$ and $\alpha: R(E)\overset{\sim}{\rightarrow}R(F)$ is an isomorphism in $\BB$. Arrows $(E,
F,\alpha)\rightarrow(E',F',\alpha')$ in $\EE\times_{\BB}\EE$ are pairs, $f: E\rightarrow E'$ and $g: F\rightarrow
F'$, such that $\alpha' R(f)=R(g) \alpha$ in $\BB$.

A fiberwise action induces a functor
\[\PP\times\EE\rightarrow \EE\times_{\BB}\EE\]
mapping $(X, E)$ to $(X\cdot E, E, b)$.

\subsection{gr-categories.} (see \cite{B1, SaRi}.)
A monoidal category $\PP$ is called \textit{grouplike} if, for each object $X$, the functor
$X\otimes(-):\PP\rightarrow\PP$ is an equivalence of categories. It follows that there is an object $X^{*}$
together with isomorphisms
\[\rho_{X}: I\rightarrow X\otimes X^{*} \mbox{ and }\lambda_{X}:X^{*}\otimes X\rightarrow I\]
satisfying the triangular identities for an adjunction (i.e.,
\[\xymatrix{X^{*}\ar[r]^-{X^{*}\otimes\rho}&X^{*}\otimes X\otimes X^{*}\ar[r]^-{\lambda\otimes X^{*}}& X^{*}}\]
and
\[\xymatrix{X\ar[r]^-{\rho\otimes X}& X\otimes X^{*}\otimes X\ar[r]^-{X\otimes \lambda}&X}\]
are both identities.) It follows also that $X\mapsto X^{*}$ can be made into a contravariant functor, for which
$\rho_{X}$ and $\lambda_{X}$ are `natural' in the appropriate sense.

A grouplike monoidal category in which each arrow $f: X\rightarrow Y$ is invertible, in other words, a grouplike
monoidal groupoid, is called a \textit{gr-category}.

\subsection{Example.}
Let $K$ be a bundle of Lie groups over a manifold $M$. Let $Bitor(M, K)$ be the monoidal category with
$K-K$-bitorsors as objects, bitorsor maps as arrows, and the tensor product defined in 3.3. Then $Bitor(M,K)$ is a
gr-category. Let $Bitor(M,T;K)$ be the monoidal category whose objects are pairs $(P,\alpha)$, where $P$ is $K-K$
torsor over $M$, and $\alpha: T\rightarrow P$ is a central section (cf. 4.3) over the submanifold $T\subset M$. A
map $(P,\alpha)\rightarrow (Q,\beta)$ in this category is a map of bitorsors $f: P\rightarrow Q$ with $f
\alpha=\beta$. The tensor product is given by $(P,\alpha)\otimes(Q,\beta)=(P\otimes Q,\alpha\otimes\beta)$. The
monoidal category thus defined is again a gr-category.

\subsection{Actions of gr-categories.}
Let $\PP$ be a gr-category, and let $R:\EE\rightarrow \BB$ be a functor. A fiberwise action of $\PP$ on
$\EE\rightarrow \BB$ is called \textit{almost principal} if the induced functor
\[\PP\times\EE\rightarrow\EE\times_{\BB}\EE\]
(cf 6.1) is an equivalence of categories. The action will be called \textit{principal} if, in addition, $R$ is
essentially surjective, in the sense that for each object $B$ in $\BB$ there is an object $E$ in $\EE$ for which
$R(E)$ is isomorphic to $B$. For such an (almost) principal action we also refer to $\EE$ as an (\textit{almost})
\textit{principal bundle} over $\BB$.

We will only use this notion when $\EE$ and $\BB$ are themselves groupoids also. In this case, there is the
following easy criterion for almost principality.
\subsection{Lemma.}
\textit{Let $R:\EE\rightarrow \BB$ be a functor between groupoids, and let $\PP$ be a gr-category with a fiberwise
action on $\EE\rightarrow \BB$. Then the action is almost principal iff the following two conditions hold:
\begin{itemize}
\item[(i)]
For any two objects $E$ and $F$ of $\EE$ and any isomorphism $\phi: R(E)\rightarrow R(F)$, there exists an object
$P$ in $\PP$ for which there exists an isomorphism $\psi: P\cdot E\rightarrow F$ in $\EE$ making the diagram
\[\xymatrix{R(P\cdot E)\ar[rr]^{b}\ar[dr]_{R\psi}&&R(E)\ar[dl]^{\phi}\\&R(F)}\]
in $\BB$ commute.
\item[(ii)]
For any $P$ in $\PP$ and $E$ in $\EE$, and any isomorphism $\psi: P\cdot E\rightarrow E$ in $\EE$ with
$R(\psi)=b$, there is a unique isomorphism $\beta: P\rightarrow I$ in $\PP$ for which the diagram
\[\xymatrix{P\cdot E\ar[rr]^{\beta\cdot E}\ar[dr]_{\psi}&&I\cdot E\ar[dl]^{u}\\&E}\]
commutes in $\EE$.
\end{itemize}}
\begin{proof}
Condition (i) is equivalent to the requirement that $\PP\times \EE\rightarrow \EE\times_{\BB}\EE$ is essentially
surjective, and (ii) is equivalent to $\PP\times \EE\rightarrow \EE\times_{\BB}\EE$ being full and faithful. We
only prove that condition (ii) implies fully faithfulness, and omit other details.

Suppose that $(P,E)$ and $(Q,F)$ are two objects in $\PP\times \EE$ and that $(\psi, \phi)$ is a map between their
images in $\EE\times_{\BB}\EE$, i.e. $\psi:P\cdot E\rightarrow Q\cdot F$ and $\phi: E\rightarrow F$ are
isomorphisms for which
\[\xymatrix{R(P\cdot E)\ar[r]^{R(\psi)}\ar[d]_{b}&R(Q\cdot F)\ar[d]^{b}\\R(E)\ar[r]^{R(\phi)}&R(F)}\]
commutes. We need to prove that there is a unique $\beta: P\rightarrow Q$ in $\PP$ such that $\psi=\beta\cdot
\phi$. First, by applying the equivalence $Q^{*}\otimes(-)$, we may assume that $Q=I$. Secondly, since $\EE$ is a
groupoid, we can compose with $\phi^{-1}: F\rightarrow E$ and assume that $E=F$ and $\phi=$identity. Finally,
composing $\psi$ with $u:I\cdot E\rightarrow E$ and observing that $b_{I,E}=R(u): R(I\cdot E)\rightarrow R(E)$, we
reduce to the problem of showing that for any isomorphism $\psi: P\cdot E\rightarrow E$ with $R(\psi)=b: R(P\cdot
E)\rightarrow R(E)$ there is a unique $\beta: P\rightarrow I$ with $\beta\cdot id_{E}=\psi$. This is condition
(ii) as stated in the lemma.
\end{proof}

Now let us go back to the specific context of the previous two sections. Let $E$ be a regular groupoid over $M$,
and let $K$ be a bundle of Lie groups over $M$. Furthermore, let $T\subset M$ be a complete transversal to the
orbit foliation of $E$. In \S4(4) we introduced the restriction functor $R: Ext_{M}(E, K)\rightarrow
Ext_{T}(E_{T}, K_{T})$. By 4.1, the gr-category $Bitor(M, K)$ of 6.3 acts on $Ext_{M}(E, K)$. The induced action
by the gr-category $Bitor(M,T;K)$ of bitorsors equipped with central sections over $T$ is a fiberwise action for
the functor $R$ (cf. Lemma 4.4). By the preceding lemma, we can now reformulate Propositions 4.5, 4.7 and 5.1 (or
5.2) together as the following theorem:

\subsection{Theorem.}
\textit{The restriction functor
\[R:Ext_{M}(E, K)\rightarrow Ext_{T}(E_{T}, K_{T})\]
is an almost principal bundle for the gr-category $Bitor(M,T;K)$. It is principal if $K$ is abelian and
$H^{2}(M,\underline{K})=0$. (A similar statement is true in the non-abelian case, cf. \S10.)}

\subsection{Remarks.}
For any discrete groupoid $\GG$ we denote by $\pi_{0}(\GG)$ the set of isomorphism class of objects. The
isomorphism class of an object $X$ in $\GG$ will be denoted $[X]$. For a gr-category $\PP$, the set $\pi_{0}(\PP)$
is a group, with multiplication induced by the tensor product.

Let $R:\EE\rightarrow \BB$ be a functor between groupoids, as earlier in this section. Let $R/ \BB$ be the comma
groupoid, with as objects triples $(E,B,\alpha)$, where $E$ is object of $\EE$, $B$ one in $\BB$, and $\alpha:
R(E)\overset{\sim}{\rightarrow}B$. Arrows $(E, B,\alpha)\rightarrow(E', B',\alpha')$ are pairs, $\phi:
E\rightarrow E'$ in $\EE$ and $u: B\rightarrow B'$ in $\BB$, such that $u \alpha=\beta R(\phi)$. There is an
evident projection functor $R/\BB\rightarrow\BB$, sending $(E, B,\alpha)$ to $B$.

If the gr-category $\PP$ acts fiberwise on $R:\EE\rightarrow\BB$, then there is an induced action of $\PP$ on the
category $R/\BB$, defined by $P\cdot(E,B,\alpha)=(P\cdot E, \xymatrix{R(P\cdot
E)\ar[r]^{b}&R(E)\ar[r]^{\alpha}&B})$.

Our definition of almost principal bundle in 6.4 is closely related to the concept of a pseudo-torsor over a site
introduced in \cite{B3}. (In fact, $R:\EE\to\BB$ is an almost principal bundle iff the fibered category $R/\BB\to
\BB$ is a pseudo-torsor over $\BB$, under the constant gr-stack $\PP \times\BB\to\BB$.)

For the proposition below we use the following terminology for group actions on sets. If a (discrete) group $G$
acts on a set $S$, and if $S\rightarrow B$ is a map which is constant on orbits, we say that $S\rightarrow B$ is
an \textit{almost principal} $G$\textit{-bundle of sets} if the action of $G$ on the fibers is free and transitive
(i.e. $G\times S\rightarrow S\times_{B}S$ is a bijection); and that $S\rightarrow B$ is a \textit{principal}
$G$\textit{-bundle} of sets if, in addition, $S\rightarrow B$ is a surjection.
\subsection{Proposition} \textit{Let $\PP$
be gr-category. If $R:\EE\rightarrow \BB$ is an (almost) principal $\PP$-bundle of groupoids, then the projection
$\pi_{0}(R/\BB)\rightarrow \pi_{0}(\BB)$ is an almost principal bundle of sets.}

The converse of this proposition is of course not true. Neither is it true that
$\pi_{0}(\EE)\rightarrow\pi_{0}(\BB)$ as an (almost) principal $\pi_{0}(\PP)$-bundle if $\EE\rightarrow\BB$ is an
(almost) principal $\PP$-bundle. This explains why we did not express the results of \S4 and Theorem 6.6 above in
terms of isomorphism classes of extensions and of bitorsors.

\subsection{Remark.}
Theorem 6.6 and the related assertions 4.5, 4.7, 5.1 and 5.2 do not involve an explicit action of $E$ on $K$.
However, statements which do are readily derived as corollaries. This is most easily explained if we assume that
$K$ is abelian. Suppose $K$ is abelian and $E$ acts on $K$ by a specific action $\mu: K\times _{M}E\rightarrow K$.
Let $Ext_{M}^{\mu}(E, K)$ be the full subgroupoid of $Ext_{M}(E, K)$ consisting of extensions $K\rightarrow
G\rightarrow E$ which induce the action this action $\mu$. Let $Ext^{\mu}_{T}(E_{T}, K_{T})$ be the similar
subgroupoid of $Ext_{T}(E_{T}, K_{T})$. Suppose a given extension $K_{T}\rightarrow B\rightarrow E_{T}$ inducing
$\mu$ is the restriction of an extension $\xymatrix{K\ar[r]^{j}&G\ar[r]^{\pi}&E}$, the latter inducing an action
$\nu$ of $E$ on $K$. Then $(K, \nu)$ is isomorphic to $(K,\mu)$ as $E$-equivariant bundles by an isomorphism
$\theta:(K, \nu)\rightarrow(K,\mu)$ which restricts to the identity on $T$, because the restriction of
$E$-equivariant bundles over $M$ to $E_{T}$-equivariant bundles over $T$ is an equivalence of categories. Thus
$\xymatrix{K\ar[r]^{j\theta}&G\ar[r]&E}$ is an extension which does induce the given action $\mu$ of $E$ on $K$,
and still restricts to $K_{T}\rightarrow B \rightarrow E_{T}$. This proves that, if $K_{T}\rightarrow B\rightarrow
E_{T}$ is up to isomorphism in the image of the restriction functor $Ext_{M}(E, K)\rightarrow Ext_{T}(E_{T},
K_{T})$, then the same is true for the functor $Ext_{M}^{\mu}(E, K)\rightarrow Ext_{T}^{\mu}(E_{T}, K_{T})$.

Next, if $K\rightarrow G\rightarrow E$ and $K\rightarrow H\rightarrow E$ are two extensions inducing the same
action $\mu$, then by Remark 4.8 the $K$-$K$ bitorsor $(P,\sigma)$ relating the two by Proposition 4.7 has the
property that $\sigma=id$. So $Ext_{M}^{\mu}(E, K)\rightarrow Ext_{T}^{\mu}(E_{T}, K_{T})$ is an almost principal
bundle for the smaller gr-category of right $K$-torsors $P$ (viewed as bitorsors $(P,id)$) with a specific section
$\alpha$ over $T$ (automatically central), and it is principal if $H^{2}(M, \underline{K})=0$.

\newpage
\part{Extensions of \'etale groupoids}
As we have seen, any regular Lie groupoid $G$ over a manifold $M$ gives rise to an extension $K_{T}\rightarrow
G_{T}\rightarrow E_{T}$ of an \'etale groupoid over a complete transversal $T$. The purpose of this part is to
introduce first a new ``\v{C}ech'' cohomology for \'etale groupoids, and then use it to describe such extensions
$K_{T}\rightarrow G_{T}\rightarrow E_{T}$. Again, we will only consider the abelian case, and defer questions
involving non-abelian cohomology to part IV.

\section{\v{C}ech cohomology of \'etale groupoids}
For the purpose of this section, let $T$ be a fixed manifold, and let $E$ be an \'etale groupoid over $T$, with
structure maps $s, t:\xymatrix{E\ar@<.5ex>[r]\ar@<-.5ex>[r]&T}$ for source and target, etc.; see \S1. Throughout
this section we also assume fixed a basis $\mathcal{U}$ for the open sets of $T$.

\subsection{Embedding categories.}
\cite{M1} Given the basis $\mathcal{U}$, the ``embedding'' category $Emb_{\mathcal{U}}(E)$ is the category whose
objects are all basis open sets $U\in \mathcal{U}$, and whose arrows $\sigma:U\to V$ are sections $\sigma: U\to E$
of the source map $s :E\to T$ such that $t \circ\sigma: U\to T$ defines an embedding $\sigma:U\to V$. Thus,
$\sigma$ can be thought of as a smooth family of arrows
\[\sigma(x):x\to t\sigma(x)\in V\qquad(\mbox{for all }x\in U).\]
For two such arrows $\sigma:U\to V$ and $\tau: V\to W$ in $Emb_{\mathcal{U}}(E)$, their composition
$\tau\circ\sigma: U\to W$ is the section defined in terms of composition in $E$,
\[(\tau\circ\sigma)(x)=\tau(t\sigma(x))\sigma(x):x\to t\sigma(x)\to t\tau(t\sigma(x)).\]
We will consider suffiently large subcategories $\KE$ of $Emb_{\mathcal{U}}(E)$: $A$ category $\KE$
\textit{approximating} the \'etale groupoid $E$ is a subcategory $\KE$ of $Emb_{\mathcal{U}}(E)$ with the
following properties
\begin{itemize}
\item[(i)] For any two basic opens $U\subseteq V$, the arrow $i:U\to V$ given by the unit section ($i(x)=1_{x}$
for $x\in U$) belongs to $\KE$. (In particular $\KE$ has the same objects as $Emb_{\mathcal{U}}(E)$).
\item[(ii)] If $\sigma: U\to V$ belongs to $\KE$, then for any $U'\subseteq U$ and $V'\subseteq V$ with
$\sigma(U')\subseteq V'$, the restriction $\sigma:U'\to V'$ also belongs to $\KE$.
\item[(iii)] For any basic open $V$ and any arrow $e:x\to y$ in $E$ with $y\in V$, there exists an arrow $\sigma:
U\to V$ in $\KE$ with $x\in U$ and $\sigma(x)=e$.
\end{itemize}

The collection of all such embedding categories $\KE$ which approximate $E$ is ``ordered'' by refinement
(inclusion) in the obvious way, and any two such embedding categories have a common refinement.

\subsection{Cohomology of small categories.} \cite{Roos} Let $\KE$ be an arbitrary small category (whose arrows
will be denoted $\sigma:U\to V$, as for the embedding categories above). Let $A$ be an abelian presheaf on $\KE$,
i.e. a contravariant functor from $\KE$ into the category of abelian groups. (We denote the effect of an arrow
$\sigma: U\to V$ by $\sigma^{*}:A(V)\to A(U)$, or by $a\mapsto a\cdot\sigma$ for $a\in A(V)$.) Then one can
construct a cochain complex $C^{\bullet}(\KE, A)$, defined in degree $n$ by
\[C^{n}(\KE,A)=\prod_{U_{0}\overset{\sigma_{1}}{\gets}U_{1}\gets\cdots\overset{\sigma_{n}}{\gets}U_{n}}A(U_{n})\]
where the product ranges over all composable strings
$U_{0}\overset{\sigma_{1}}{\gets}U_{1}\gets\cdots\overset{\sigma_{n}}{\gets}U_{n}$ of arrows in $\KE$. The
coboundary $d:C^{n-1}(\KE,A)\to C^{n}(\KE,A)$ is defined, for $c\in C^{n}(\KE,A)$, by $d
c=\sum_{i=0}^{n}(-1)^{i}d_{i}(c)$, where
\[d_{i}(c)(U_{0}\overset{\sigma_{1}}{\gets}\cdots\overset{\sigma_{n}}{\gets}{U_{n}})
=\left\{\begin{array}{lr}c(U_{1}\overset{\sigma_{2}}{\gets}\cdots\gets U_{n})&\mbox{if }i=0,\\
c(U_{0}\overset{\sigma_{1}}{\gets}\cdots U_{i-1}\overset{\sigma_{i}\sigma_{i+1}}{\gets}U_{i+1}\cdots\gets U_{n})&\mbox{if }0<i<n,\\
c(U_{0}\overset{\sigma_{1}}{\gets}\cdots\overset{\sigma_{n-1}}{\gets}U_{n-1})\cdot\sigma_{n}&\mbox{if }i=n.
\end{array}\right.\]
The cohomology of this complex is written $H^{*}(\KE,A)$, and referred to as the cohomology of $\KE$ with
coefficients in $A$.

\subsection{Cohomology of embedding categories.}
Now let $E$ be an \'etale groupoid over $T$, as above, and let $A$ be a sheaf on $T$ with a right action by $E$.
There is a standard sheaf cohomology $H^{*}(E, A)$, first considered by Haefliger. It agrees with the cohomology
of the classifying space $BE$, in the sense that
\begin{equation}\label{coh:eq1}
H^{*}(E,A)=H^{*}(BE,\tilde{A})
\end{equation}
where $\tilde{A}$ is the sheaf on $BE$ induced by $A$ in a canonical way. See \cite{M2} for a proof of
(\ref{coh:eq1}).

Such an $E$-equivariant sheaf $A$ on $T$ defines a presheaf $\Gamma(A)$ on $\KE$, defined on objects by
\[\Gamma(A)(U)=\Gamma(U,A),\]
and arrows $\sigma: U\to V$ of $\KE$ by
\begin{equation}\label{coh:eq2}
(a\cdot\sigma)(x)=a(t\sigma(x))\cdot\sigma(x),
\end{equation}
where $a$ is a section of $A$ over $V$, and $x$ is a point in $U$; the dot on the right of (\ref{coh:eq2}) refers
to the action of $E$ on $A$. We recall the following result from \cite{M1, CM}. (In fact, the following statement
is slightly more general than the one given in \cite{CM}, but the same proof applies.)

\subsection{Theorem.}
\textit{Let ${\KE}$ be an embedding category approximating the \'etale groupoid $E$, and let $A$ be an
$E$-equivariant sheaf of abelian groups on $T$. Suppose that each object $U$ of $\KE$ has the property that
$H^{i}(U, A)=0$ for $i>0$. Then
\[H^{*}(E, A)=H^{*}(\KE, \Gamma(A)).\]}

\subsection{\v{C}ech cohomology.}
We will think of an embedding category $\KE$ approximating $E$ as a cover of $E$, and $H^{*}(\mathcal{E},\Gamma
A)$ as the \v{C}ech cohomology of $E$ for this cover, so we also write
\[\check{H}^{*}(\KE, A)=H^{*}(\KE,\Gamma A).\]
If $\KE'$ is an embedding which refines $\KE$, there is an evident restriction map $\check{H}^{*}(\KE, A)\to
\check{H}^{*}(\KE', A)$. The \v{C}ech cohomology of $E$ with coefficients in $A$ is now defined as the colimit
over all embedding categories approximating $E$,
\[\check{H}^{*}(E, A)=colim_{\KE}\check{H}^{*}(\KE, A).\]
Thus, as a consequence of Theorem 7.4, we find that this \v{C}ech cohomology agrees with sheaf cohomology under
suitable conditions; for example
\[\check{H}^{*}(E, A)\simeq H^{*}(E, A),\]
if $A$ is an $E$-equivariant sheaf which is locally constant as a sheaf on $T$.

\subsection{Remark.}
For an arbitrary $E$-equivariant abelian sheaf $A$ on $T$, there is always a map $\check{H}^{*}(E, A)\to H^{*}(E,
A)$ from \v{C}ech cohomology to the sheaf cohomology. The \v{C}ech cohomology introduced here is a ``finer''
approximation of sheaf cohomology than the one introduced by Haefliger in \cite{H1}, in the sense that the
canonical map $\check{H}_{Haefliger}^{*}(E, A)\to H^{*}(E,A)$ factors through the previous map. Like Haefliger's,
our \v{C}ech cohomology is invariant under Morita equivalence of \'etale groupoids (cf \S1).

\subsection{Relative \v{C}ech cohomology.}
If $\KE$ is a small category and $\mathcal{D}\subset \KE$ is a subcategory, we write
\[C^{n}(\KE,\mathcal{D};A)\]
for the group of $n$-cochains on $\KE$ with values in the presheaf $A$ which vanish on $\mathcal{D}$. These
cochains form a subcomplex of $C^{*}(\KE, A)$ which fits into a short exact sequence
\[0 \to C^{*}(\KE,\mathcal{D};A)\to C^{*}(\KE,A)\to C^{*}(\mathcal{D},A)\to 0.\]
In cohomology this gives a long exact sequence
\begin{equation}\label{Cech:eq1}
\cdots\to H^{n}(\KE,\mathcal{D};A)\to H^{n}(\KE,A)\to H^{n}(\mathcal{D},A)\to H^{n+1}(\KE,\mathcal{D};A)\to\cdots,
\end{equation}
as usual.

For an embedding category $\KE$ of the \'etale groupoid $\xymatrix{E\ar@<.5ex>[r]\ar@<-.5ex>[r]&T}$, we write
$\KE_{u}$ for the subcategory with the same objects as $\KE$, and with only those arrows $\sigma: U\to V$ for
which $\sigma$ is a unit section (i.e. $\sigma(x)=1_{x}$ for all $x\in U$, and $U\subset V$). For an
$E$-equivariant sheaf $A$, the presheaf $\Gamma A$ on $\KE_{u}$, and the cohomology $H^{n}(\KE_{u},\Gamma A)$ -
which we also write as $\check{H}^{n}(\KE_{u},A)$ - is some sort of \v{C}ech cohomology of the manifold $T$. If
the sheaf $A$ satisfies the conditions of Theorem 7.4 then this theorem applied to $\KE_{u}$ gives an isomorphism
$H^{*}(T,A)=\check{H}^{*}(\KE_{u},A)$. For general $A$, we write
\[\check{H}^{*}(T,A)=colim_{\KE}\check{H}^{*}(\KE_{u},A).\]
Let us write $\check{H}^{n}(E,T;A)$ for the relative \v{C}ech cohomology $colim_{\KE}(\KE,\KE_{u};A)$. Thus, by
passing to the colimit over approximating embedding categories, the long exact sequence (\ref{Cech:eq1}) with
$\KE_{u}$ for $\mathcal{D}$ yields a long exact sequence
\[\cdots\to \check{H}^{n}(E,T;A)\to \check{H}^{n}(E,A)\to \check{H}^{n}(T,A)\to\check{H}^{n+1}(E,T;A)\to\cdots\]
Intuitively speaking, one should think of the inclusion $T\hookrightarrow E$ of the units into the \'etale
groupoid as a quotient map. It is a model for the ``homotopy quotient'' $T\to BE$. The long exact sequence above
is the one associated to this quotient map.

\subsection{Remark.}
Suppose the sheaf $A$ satisfies the conditions of Theorem 7.4, so that the l.e.s. takes the form
\[\cdots\to  H^{n-1}(T,A)\to H^{n}(E,T,A)\to H^{n}(E,A)\to H^{n}(T,A)\to\cdots\]
Thus, if $H^{n}(T,A)=0$ for $n>0$ then $H^{n}(E,T,A)=H^{n}(E,A)$ for $n\geq 2$. Let $\mathcal{U}$ be an open cover
of $T$ by contractible open sets, and let $T'=\coprod_{U\in\mathcal{U}}U$ be the disjoint sum. Then $E$ pulls back
to an \'etale groupoid $E'$ over $T'$, which is Morita equivalent to $E$ and has the property that $H^{n}(T',A)=0$
for $n>0$. Thus $H^{n}(E',T',A)=H^{n}(E',A)=H^{n}(E,A)$ for $n\geq 2$, the latter identity by Morita invariance.
We will use this observation in 8.7 below, to relate our classification theorem tot that of Kumjian et al.

\section{Extensions and cohomology}
In this section we show how extensions of \'etale groupoids are classified by degree 2 cohomology classes in the
\v{C}ech cohomology described in the previous section. We deal with the abelian case first. In Section 12 below we
will show that a completely analogous classification holds for an appropriately defined non-abelian cohomology.

As in the previous section, we fix an \'etale groupoid $\xymatrix{E\ar@<.5ex>[r]\ar@<-.5ex>[r]&T}$ over a manifold
$T$. Let $K$ be a bundle of abelian Lie groups over $T$, equipped with a right action $\mu:K\times_{T}E\to K$,
which we simply denote for $e:x\to y$ and $k\in K_{y}$ by $\mu(k,e)=k\cdot e=k e$. The sheaf of smooth sections of
$K$ will be denoted $\underline{K}$. It is an $E$-equivariant sheaf on $T$.

\subsection{The group $\pi_{0} Ext^{\mu}_{T}(E,K)$.}
Recall from 6.8 that $Ext^{\mu}_{T}(E,K)$ denotes the full subgroupoid of the groupoid $Ext_{T}(E,K)$ of
extensions $K\to G\to E$ over $T$, consisting of only those extensions which induce the given action $\mu$. Since
$K$ is assumed abelian, this subcategory $Ext^{\mu}_{T}(K,E)$ has the structure of a gr-category (6.2), the tensor
product being given by the ``Baer sum'' of extensions. Explicitly, for two extensions
$K\overset{i}{\to}G\overset{\pi}{\to}E$ and $K\overset{j}{\to}H\overset{\rho}{\to}E$, this Baer sum is the
extension
\begin{equation}\label{8:1}
K\to G\otimes_{K}H\to E
\end{equation}
where $G\otimes_{K}H$ is the Lie groupoid whose arrows $x\to y$ are equivalence classes $g\otimes h$ of pairs
$x\overset{g}{\to}y$, $y\overset{h}{\to}x$ with $g$ and $h$ arrows in $G$ and $H$, respectively, such that
$\rho(h)\pi(g)=1_{x}.$ The equivalence relation on such pairs is given by the identification $g k\otimes
h=g\otimes k h$ --- or more precisely, $g i(k)\otimes h=g\otimes j(k) h$, for any $k\in K_{x}$. The composition in
the groupoid $G\otimes_{K}H$ is defined as
\[(g \otimes h)(g'\otimes h')=g g'\otimes h' h.\]
(This is well defined on equivalence classes because $G$ and $H$ induce the same action $\mu$ of $E$ on $K$.) The
maps $K\to G\otimes_{K}H$ and $G\otimes_{K}H\to E$ in the extension (\ref{8:1}) send $k\in K_{x}$ to $i(k)\otimes
1_{x}=1_{x}\otimes i(k)$, and $g\otimes h$ to $\pi(g)=\pi(h^{-1})$. The unit for the Baer sum is the semidirect
product groupoid $K\rtimes E$. Its arrows are pairs $(k,e)$ with $k\in K_{x}$ and $e:x\to y$ in $E$. Its
composition is defined by the usual formula involving the action: for $(l,f):y\to z$ in $K\rtimes E$,
\[(l,f)(k,e)=((l\cdot e)k, f e).\]
The inverse of an extension $K\overset{i}{\to}G\overset{\pi}{\to}E$ is the extension
$K\overset{i^{op}}{\to}G^{op}\overset{\pi^{op}}{\to}E$, where $G^{op}$ is the opposite groupoid of $G$ (arrows
$g^{op}:x\to y$ in $G^{op}$ are arrows $g:y\to x$ in $G$), while $i^{op}(k)=i(k^{-1})^{op}$ and
$\pi^{op}(g^{op})=\pi(g^{-1})$. (This extension $G^{op}$ indeed induces the same action $\mu$.) Observe that the
B\"ar sum is in fact symmetric, by the isomorphism
\[G\otimes_{K}H\overset{\sim}{\to}H\otimes_{K}G\]
mapping $g\otimes h$ to $h^{-1}\otimes g^{-1}$. Thus $Ext_{T}^{\mu}(E,K)$ is in fact a symmetric gr-category, i.e.
a \textit{Picard category}. It follows that the set of isomorphism classes of objects in this category,
\[\pi_{0}Ext^{\mu}_{T}(E,K)\]
has the structure of an abelian group.

\subsection{Theorem.}
\textit{There is a natural isomorphism of groups
\[\check{H}^{2}(E,T;\underline{K})\simeq\pi_{0}Ext^{\mu}_{T}(E,K).\]}
We prove this theorem by explicitly exhibiting the construction of a \v{C}ech cocycle from an extension (8.3) and
of an extension from a cocycle (8.5). It will be obvious (8.6) that these constructions are mutually inverse and
give the stated isomorphism.

\subsection{From extensions to cocycles.}
Suppose we are given an extension over $T$,
\[K\overset{j}{\to}G\overset{\pi}{\to}E\]
An embedding category is said to be \textit{adapted} to this extension if, for every arrow $\sigma:U\to V$ in
$\KE$, there exists a section $\tilde{\sigma}:U\to G$ of the source map $G\to T$, such that
$\pi\circ\tilde{\sigma}=\sigma$. In other words, for each $x\in U$ the arrow $\sigma(x):x\to t\sigma(x)=y$ is
lifted to an arrow $\tilde{\sigma}(x):x\to y$ in $G$. Since $\pi$ is a surjective submersion, there exist
arbitrary fine embedding categories $\KE$ which approximate $E$ and are adapted to the given extension. Let us fix
one such embedding category $\KE$. Also, let us fix for each $\sigma:U\to V$ in this $\KE$ a lifting
$\tilde{\sigma}$ as above. In case $\sigma:U\to V$ is a unit section of $E$ (in the sense that $\sigma(x)=1_{x}$
for all $x\in U$), we agree to choose $\tilde{\sigma}$ to be the corresponding unit section of $G$. Now define a
\v{C}ech cocycle $c$ by
\[c(U_{0}\overset{\sigma_{1}}{\gets}U_{1}\overset{\sigma_{2}}{\gets}U_{2})(x)=(\widetilde{\sigma_{1}\sigma_{2}})(x)^{-1}
\tilde{\sigma_{1}}(y)\tilde{\sigma_{2}}(x)\] where $x\in U_{2}$ and $y=t \sigma_{2}(x)$. Then
$c(\sigma_{1},\sigma_{2})(x)\in K_{x}$, so $c(\sigma_{1},\sigma_{2})\in\Gamma(U_{2},\underline{K})$ is a section
of the sheaf $\underline{K}$.

To see that $c$ is indeed a cocycle, consider the sequence
\[U_{0}\overset{\sigma_{1}}{\gets}U_{1}\overset{\sigma_{2}}{\gets}U_{2}\overset{\sigma_{3}}{\gets}U_{3}.\]
For a point $x\in U_{3}$, let us write $y=t\sigma_{3}(x)$ and $z=t\sigma_{2}(y)$. Let us also write $\sigma_{1
2}=\sigma_{1}\circ\sigma_{2}$, etc. Then
\[\begin{array}{ll}
c(\sigma_{1},\sigma_{2 3})(x) c(\sigma_{2},\sigma_{3})(x)&=\tilde{\sigma}_{1 2 3}(x)^{-1}\tilde{\sigma}_{1}(z)
\tilde{\sigma}_{2 3}(x)\tilde{\sigma}_{2 3}^{-1}(x)\tilde{\sigma}_{2}(y)\tilde{\sigma}_{3}(x)\\
&=\tilde{\sigma}_{1 2 3}(x)^{-1}\tilde{\sigma}_{1}(z)\tilde{\sigma}_{2}(y)\tilde{\sigma}_{3}(x),
\end{array}\]
while also
\[\begin{array}{l}
c(\sigma_{1 2},\sigma_{3})(x)(c(\sigma_{1},\sigma_{2})(y)\cdot\sigma_{3}(x))=\\
\qquad=[\tilde{\sigma}_{1 2 3}(x)^{-1}\tilde{\sigma}_{1 2}(y)\tilde{\sigma}_{3}(x)]\tilde{\sigma}_{3}(x)^{-1}
[\tilde{\sigma}_{1 2}(y)^{-1}\tilde{\sigma}_{1}(z)\tilde{\sigma}_{2}(y)]\tilde{\sigma}_{3}(x)\\
\qquad=\tilde{\sigma}_{1 2 3}(x)^{-1}\tilde{\sigma}_{1}(z)\tilde{\sigma}_{2}(y)\tilde{\sigma}_{3}(x).
\end{array}\]
Thus $c(\sigma_{1},\sigma_{2 3}) c(\sigma_{2},\sigma_{3})=c(\sigma_{1
2},\sigma_{3})(c(\sigma_{1},\sigma_{2})\cdot\sigma_{3})$, whence $(d c)(\sigma_{1},\sigma_{2},\sigma_{3})=0$. This
proves that $c$ is a cocycle. Also notice that if $\sigma_{1},\sigma_{2}$ are both unit sections, then by our
choice so are $\tilde{\sigma}_{1}$ and $\tilde{\sigma}_{2}$, so $c(\sigma_{1},\sigma_{2})(x)=1_{x}$ for all $x\in
U_{2}$. Thus $c$ is in fact a relative cocycle (7.7), and hence define a cohomology class
\[ [c]\in \check{H}^{2}(\KE,\KE_{u},\underline{K}).\]

This cohomology class is independent of the chosen liftings $\tilde{\sigma}$ for each $\sigma$. Indeed, suppose
$\bar{\sigma}$ is another such choice for each arrow $\sigma$ in $\KE$, and let $\bar{c}$ be the corresponding
cocycle. We can then define a cochain $f\in C^{1}(\KE,\KE_{u},\underline{K})$ by
\[f(U_{0}\overset{\sigma}{\gets}U_{1})(y)=\bar{\sigma}(y)^{-1}\tilde{\sigma}(y)\in K_{y} \qquad
(\mbox{for }y\in U_{1}).\] Then, for $U_{0}\overset{\sigma_{1}}{\gets}U_{1}\overset{\sigma_{2}}{\gets}U_{2}$ and
$x\in U_{2}$ with $y=t\sigma_{2}(x)$, we have
\[\begin{array}{ll}
f(\sigma_{1 2})(x)c(\sigma_{1},\sigma_{2})(x)&=\bar{\sigma}_{1 2}(x)^{-1}\tilde{\sigma}_{1 2}(x)
\tilde{\sigma}_{1 2}(x)^{-1}\tilde{\sigma}_{1}(y)\tilde{\sigma}_{2}(x)\\
&=\bar{\sigma}_{1 2}(x)^{-1}\tilde{\sigma}_{1}(y)\tilde{\sigma}_{2}(x),
\end{array}\]
while
\[\begin{array}{l}
\bar{c}(\sigma_{1},\sigma_{2})(x)(f(\sigma_{1})(y)\cdot\sigma_{2}(x))f(\sigma_{2})(x)=\\
\qquad =[\bar{\sigma}_{1 2}(x)^{-1}\bar{\sigma}_{1}(y)\bar{\sigma}_{2}(x)][\bar{\sigma}_{2}(x)^{-1}
\bar{\sigma}_{1}(y)^{-1}\tilde{\sigma}_{1}(y)\bar{\sigma}_{2}(x)][\bar{\sigma}_{2}(x)^{-1}\tilde{\sigma}_{2}(x)]\\
\qquad =\bar{\sigma}_{1 2}(x)^{-1}\tilde{\sigma}_{1}(y)\tilde{\sigma}_{2}(x),
\end{array}\]
also. Thus
\[f(\sigma_{1 2}) c(\sigma_{1},\sigma_{2})=\bar{c}(\sigma_{1},\sigma_{2})(f(\sigma_{1})\cdot\sigma_{2})f(\sigma_{2}),\]
which expresses in the abelian case that $c$ and $\bar{c}$ differ by the coboundary $d f$.

By mapping into the colimit over the embedding categories $\KE$ adapted to $G$, we obtain in this way a
well-defined class
\[ [c]\in \check{H}^{2}(E,T;\underline{K})\]
associated to the extension $K\to G\to E$.

\subsection{Normalization.}
It is possible to choose a lifting  $\tilde{\sigma}$ of each arrow $\sigma:U\to V$ in $\KE$ in a special way. We
already  agreed that if $\sigma(x)=1_{x}$ for all $x\in U$ (i.e. $\sigma$ is a unit section) then the same is true
for $\tilde{\sigma}$. Next, each arrow $\sigma:U\to V$ in $\KE$ factors as an isomorphism $\sigma_{0}$ followed by
a unit section $i$,
\[U\overset{\sigma_{0}}{\to}t \sigma(U)\overset{i}{\hookrightarrow}V,\]
by $\sigma_{0}(x)=\sigma(x)$ and $i(y)=1_{y}$. We can fix a lifting $\tilde{\sigma}_{0}$ of
$\sigma_{0}:U\overset{\sim}{\to}t \sigma(U)$, and define $\tilde{\sigma}$ by
$\tilde{\sigma}(x)=\tilde{\sigma}_{0}(x)$. In this way, we obtain a cocycle with the property that
\[c(U_{0}\overset{\sigma_{1}}{\gets}U_{1}\overset{\sigma_{2}}{\gets}U_{2})=1\]
whenever $\sigma_{1}$ is a unit section. We call this cocycle \textit{normal}. In general, any cocycle $\xi\in
C^{2}(\KE,\KE_{u},\underline{K})$ can be normalized. Indeed, consider such a relative 2-cocycle $\xi$, so that
\[\xi(j_{1},j_{2})=1\]
whenever $U_{0}\overset{j_{1}}{\gets}U_{1}\overset{j_{2}}{\gets}U_{2}$ are two unit sections. Let $g$ be the
1-cochain defined by
\[g(\sigma)=\xi(i,\sigma_{0}),\]
where $\sigma=i \sigma_{0}$ with $\sigma_{0}$ an isomorphism and $i$ a unit section. Then $\xi$ is cohomologous to
the 2-cocycle $\hat{\xi}$ defined on $U_{0}\overset{\tau}{\gets}U_{1}\overset{\sigma}{\gets}U_{2}$ by
\[\hat{\xi}(\tau,\sigma)=g(\tau\sigma)^{-1}\xi(\tau,\sigma)(g(\tau)\cdot\sigma)g(\sigma).\]
Now $\hat{\xi}$ is normal, because if $\tau=j$ is a unit section, and $\sigma=i \sigma_{0}$ as above, the cocycle
condition for $\xi$ applied to $\overset{j}{\gets}\overset{i}{\gets}\overset{\sigma_{0}}{\gets}$ gives
\[\xi(j,\sigma)\xi(i,\sigma_{0})=\xi(j i,\sigma_{0})(\xi(j,i)\cdot\sigma_{0}).\]
Since $\xi$ vanishes on $\KE_{u}$, we can rewrite this as
\begin{equation}\label{8:2}
\xi(j,\sigma)\xi(i,\sigma_{0})=\xi(j i,\sigma_{0}).
\end{equation}
But then
\[\begin{array}{ll}\hat{\xi}(j,\sigma)&=g(j \sigma)^{-1}\xi(j,\sigma)(g(j)\cdot\sigma)g(\sigma)\\
&=g(j \sigma)^{-1}\xi(j,\sigma)g(\sigma)\qquad(\mbox{because }g(j)=1)\\
&=\xi(j i,\sigma_{0})^{-1}\xi(j,\sigma)\xi(i,\sigma_{0})\\
&=1\qquad\mbox{by (\ref{8:2})}.
\end{array}\]
For later calculations, we observe that normal cocycles $\zeta$ have the property that
\begin{equation}\label{rr3}
\zeta(j \sigma,\tau)=\zeta(\sigma,\tau)
\end{equation}
whenever $j$ is a unit section, and
\begin{equation}\label{rr4}\zeta(\sigma,1)=1\end{equation}
where $1$ is a unit isomorphism, i.e. an identity arrow in $\KE$. Notice also that the cocycle condition implies
that
\begin{equation}\label{rr5}\zeta(\sigma,\sigma^{-1})\cdot\sigma=\zeta(\sigma^{-1},\sigma)\end{equation}
for any cocycle $\zeta$.

\subsection{From cocycles to extensions.}
Assume now we are given a cohomology class $[\xi]\in\check{H}^{2}(E,T;\underline{K})$, which we may take to be
represented by a $\textit{normal}$ relative cocycle
\[\xi\in C^{2}(\KE,\KE_{u};\underline{K})\]
for some embedding category $\KE$ which approximates $E$.

We now construct a groupoid $G(\xi)$ which fits into an extension
\begin{equation}\label{rr6}
K\overset{\nu}{\to}G(\xi)\overset{\pi}{\to}E
\end{equation}
of groupoids over $T$. For two points $x,y\in T$, arrows $x\to y$ in $G(\xi)$ are equivalence classes of pairs
$(\sigma,k)$ where $\sigma:U\overset{\sim}{\to} V$ is an isomorphism in the embedding category whose domain is a
nbd $U$ of $x$, such that $t\sigma(x)=y$, and where $k\in K_{x}$. For a smaller nbd $U'\subset U$ of $x$, the pair
$(\sigma,k)$ is identified with the pair $(\sigma',\xi(\sigma,i)(x) k)$, where $\sigma':U'\overset{\sim}{\to} V'$
is the restriction of $\sigma$, $V'=\sigma(U')$, and $i:U'\to U$ is the arrow in $\KE$ given by the unit section.
The equivalence class of $(\sigma,k)$ is denoted
\[[\sigma,k].\]
The maps $\nu$ and $\pi$ in (20) are defined by
\[\begin{array}{l}
\nu(k)=[1_{U},k],\\
\pi[\sigma,k]=\sigma(x).
\end{array}\]
(In this definition of $\nu$, $U$ is any nbd of $x$ and $1_{U}$ is the identity at $U$ in $\KE$; the class
$[1_{U},k]$ is independent of the choice of $U$, because $\xi(1_{U},i)=1$ is the unit section of $K$.)

Before we go into the groupoid structure of $G(\xi)$, let us first observe that $G(\xi)$ is a smooth manifold. If
$U\subset T$ is an open set and $\sigma: U\overset{\sim}{\to}V$ is an isomorphism in $\KE$, then $\{\sigma(x)|x\in
U\}$ is an open subset of $E$, diffeomorphic to $U$. Write $K_{U}=\bigcup_{x\in U}K_{x}$. Then $K_{U}$ is an open
subset of $K$. There is a canonical map $K_{U}\to G(\xi)$ which sends $k\in K_{x}$ to $[\sigma,k]:x\to
t\sigma(x)$. This map is injective, and we define the smooth structure on $G(\xi)$ by declaring $K_{U}\to G(\xi)$
to be an open embedding. In other words, whereas $E$ is the union of copies of open sets $U\subset T$, glued
together via the sections $\sigma$, the space $G(\xi)$ is a union of copies of $K_{U}$, glued together in the same
pattern.

To gain some feeling fore the space $G(\xi)$, it is useful also to observe that, for a point $x\in T$, the space
$G(\xi)(x,-)$ of arrows out of $x$ is simply the product $E(x,-)\times K_{x}$. Indeed, the space $E(x,-)$ is
discrete because $E$ is \'etale. So by choosing a section $\sigma_{e}$ through each $e\in E(x,-)$, we can define a
diffeomorphism $E(x,-)\times K_{x}\to G(\xi)(x,-)$ by mapping $(e,k)$ to $[\sigma_{e},k]$.

Composition in $G(\xi)$ is defined as follows. For two arrows $[\sigma,k]:x\to y$ and $[\tau,l]:y\to z$, we may
assume that $\sigma:U\overset{\sim}{\to}V$ and $\tau:V\overset{\sim}{\to}W$, and then define
\begin{equation}\label{rr7}
[\tau,l][\sigma,k]=[\tau\sigma,\xi(\tau,\sigma)(x)(l\cdot\sigma(x)k]
\end{equation}
where $l\cdot\sigma(x)$ refers to the action of $E$ on $K$. We need to prove, first of all, that this is
well-defined on equivalence classes. To this end, consider smaller nbds $U\subset U'$, $V\subset V'$ and
$W'\subset W$ for which $\sigma$ and $\tau$ restrict to $\sigma':U'\overset{\sim}{\to} V'$ and
$\tau':V'\overset{\sim}{\to} W'$. Also, write $U'\overset{i}{\to}U$, $V'\overset{j}{\to}V$ and
$W'\overset{h}{\to}W$ for the inclusions as arrows in the embedding category:
\[\xymatrix{W'\ar[d]_{h}&V'\ar[l]^{\tau'}\ar[d]_{j}&U'\ar[l]^{\sigma'}\ar[d]^{i}\\
W&V\ar[l]^{\tau}&U\ar[l]^{\sigma}}\] Then
\[[\sigma,k]=[\sigma',\xi(\sigma,i)(x) k]\mbox{ and }[\tau,l]=[\tau',\xi(\tau,j)(y) l]\]
So, for well-definedness of the composition (21), we need to show that
\[\begin{array}{l}\xi(\tau \sigma,i)(x)\xi(\tau,\sigma)(x)(l\cdot\sigma(x))k=\\
\qquad =\xi(\tau',\sigma')(x)[(\xi(\tau,j)(y) l)\cdot\sigma'(x)]\xi(\sigma,i)(x) k.\end{array}\] Since
$\sigma(x)=\sigma'(x)$ and $K$ is assumed abelian, it suffices to show that
\begin{equation}\label{rr8}
\xi(\tau\sigma,i)(\xi(\tau,\sigma)\cdot i)=\xi(\tau',\sigma')(\xi(\tau,j)\cdot\sigma')\xi(\sigma,i)
\end{equation}
as sections of $K$ on $U'$. Now notice that, by the cocycle condition for
$\overset{\tau}{\gets}\overset{j}{\gets}\overset{\sigma'}{\gets}$, we have
\[\xi(\tau,j\sigma')\xi(j,\sigma')=\xi(\tau j,\sigma')(\xi(\tau, j)\cdot\sigma').\]
By normality of $\xi$, this reduces to
\[\xi(\tau,j)\cdot\sigma'=\xi(\tau j,\sigma')^{-1}(\xi(\tau,j \sigma').\]
Now $\xi(\tau j,\sigma')=\xi(h \tau',\sigma')=\xi(\tau',\sigma')$ (the latter by (17)), while
$\xi(\tau,j\sigma')=\xi(\tau,\sigma i)$. so the r.h.s. of (22) equals
\[\xi(\tau,\sigma i)\xi(\sigma,i).\]
By the cocycle condition for $\overset{\tau}{\gets}\overset{\sigma}{\gets}\overset{i}{\gets}$, this last
expression equals $\xi(\tau\sigma,i)(\xi(\tau,\sigma)\cdot i)$, which is exactly the l.h.s. of (22). This proves
composition is well-defined on equivalence classes.

A straightforward calculation using the cocycle condition on $\xi$ now shows that composition of $G(\xi)$ thus
defined is associative, and that each arrow $[\sigma, k]:x\to y$ has an inverse, namely $[\sigma^{-1},l]:y\to x$
where $l=\xi(\sigma,\sigma^{-1})^{-1}(k^{-1}\cdot\sigma^{-1})$. (use (19) for this). We omit the details. Finally,
composition thus defined makes $\nu:K\to G(\xi)$ and $\pi:G(\xi)\to E$ into groupoid homomorphisms, and provides
the exact sequence (20).

Finally, we note that, up to isomorphism, the construction of the groupoids $G(\xi)$ depends only on the
cohomology class of the cocycle $\xi$. Indeed if $\zeta$ is another normal cocycle, cohomologous to $\xi$, we can
write
\[\xi(\tau,\sigma)=g(\tau\sigma)^{-1}\xi(\tau,\sigma)(g(\tau)\cdot\sigma)g(\sigma)\]
for a 1-cochain $g$, which we may assume to be normal in the sense that $g(1_{U})=1$ for any identity arrow
$1_{U}:U\to U$ in $\KE$. But then one easily verifies that the map
\[G(\xi)\to G(\zeta)\]
sending $[\sigma,k]$ to $[\sigma,g(\sigma)^{-1}k]$ is a groupoid homomorphism, which makes $K\to G(\xi)\to E$ and
$K\to G(\zeta)\to E$ isomorphic as extensions.

\subsection{The constructions are mutually inverse.}
To conclude the proof of Theorem 8.2, we briefly indicate why the constructions of 8.3 and 8.5 are mutually
inverse. One way round, suppose we are given an extension $K\overset{i}{\to}G\overset{\pi}{\to}E$. For an
embedding category $\KE$ adapted to this extension, we chose liftings $\tilde{\sigma}$ of each arrow $\sigma:U\to
V$; we can assume that this was done in a ``normal way'', as in 8.4, so that $\widetilde{j\sigma}=\tilde{\sigma}$
if $j$ is a unit section. This gave rise to a normal cocycle $\xi$, denoted informally by
$\xi(\tau,\sigma)=\widetilde{\tau\sigma}^{-1}\tilde{\tau}\tilde{\sigma}$. From this cocycle $\xi$, we subsequently
constructed a groupoid $G(\xi)$ in 8.5. To see that $G$ is isomorphic to $G(\xi)$ (as extensions), define a map
\[\phi:G\to G(\xi)\]
by setting, for an arrow $g:x\to y$ in $G$,
\[\phi(g)=[\sigma,\tilde{\sigma}(x)^{-1}g],\]
where $\sigma:U\overset{\sim}{\to}V$ is an isomorphism in $\KE$ with $\sigma(x)=\pi(g):x\to y$ in $E$. Such a
$\sigma$ exists because $\KE$ approximates $E$, but $\sigma$ is not unique. However, any two choices for $\sigma$
must agree on a sufficiently small nbhd of $x$, and from this it follows that $\phi(g)$ is well-defined and
independent of the choice of $\sigma$. Indeed, if $U'\subset U$ and $\sigma':U'\overset{\sim}{\to}V'$ denotes the
restriction of $\sigma$ to a smaller nbhd $U'$ of $x$, then if we write $i:U'\to U$ and $j:V'\to V$ for the
inclusions in $\KE$, we have in $G(\xi)$ the identity
\[\begin{array}{ll}
[\sigma,\tilde{\sigma(x)}^{-1}g]&=[\sigma',\xi(\sigma,i)(x)\tilde{\sigma}(x)^{-1}g]\\
&=[\sigma',\tilde{\sigma}'(x)^{-1}g],
\end{array}\]
the latter because
\[\begin{array}{ll}
\xi(\sigma,i)(x)&=(\widetilde{\sigma i})(x)^{-1}\tilde{\sigma}(x)i(x)\\
&=(\widetilde{j\sigma'})(x)^{-1}\tilde{\sigma}(x)\\
&=\tilde{\sigma}'(x)^{-1}\tilde{\sigma}(x),
\end{array}\]
by normality of $\xi$.

One now easily checks that this map $\phi: G\to G(\xi)$ preserves composition, and provides the required
isomorphism of extensions.

The other way round, suppose we are given a normal cocycle defined on an approximating embedding category $\KE$.
Form the groupoid $G(\xi)$ as in 8.5. The same category $\KE$ will then be adapted to the extension $K\to
G(\xi)\to E$, and as a lifting $\tilde{\sigma}$ of $\sigma:U\overset{\sim}{\to}V$ in $\KE$ we may take
\[\tilde{\sigma}(x)=[\sigma,1_{x}]:x\to y.\]
For an arbitrary $\sigma:U\to W$ in $\KE$ we factor $\sigma$ as an isomorphism $\sigma_{0}:U\overset{\sim}{\to}V$
followed by an inclusion $i:V\hookrightarrow W$, and define $\tilde{\sigma}(x)=[\sigma_{0},1_{x}]$. Then this
choice of liftings defines a normal cocycle $\zeta$,
\[\zeta(\tau,\sigma)=(\widetilde{\tau\sigma})^{-1}\tilde{\tau}\tilde{\sigma},\]
or more precisely,
\begin{equation}\label{rr9}
\zeta(\tau,\sigma)(x)=[\tau_{1}\sigma_{0},1_{x}]^{-1}[\tau_{0},1_{y}][\sigma_{0},1_{x}],
\end{equation}
where $\sigma=i\sigma_{0}$, $\tau=k\tau_{0}$ and $\tau_{0}i=j\tau_{1}$ so that $\tau\sigma=(k
j)(\tau_{1}\sigma_{0})$; here $\sigma_{0}$, $\tau_{0}$ and $\tau_{1}$ are isomorphisms and $i$, $j$ and $k$ are
inclusions (unit sections of $\KE$), as in the diagram
\[\xymatrix{U\ar[d]^{\sim}_{\sigma_{0}}\ar[r]^{\sigma}&W\ar@{=}[d]\ar[r]^{\tau}&N\ar@{=}[d]\\
V\ar[d]^{\sim}_{\tau_{1}}\ar[r]^{i}&W\ar[d]_{\tau_{0}}^{\sim}\ar[r]^{\tau}& N\ar@{=}[d]\\
0\ar[r]^{j}&M \ar[r]^{k}&N.}\]

 We claim that the cocycle $\zeta$ defined by (23) is in fact identical to $\xi$.
Indeed, $[\tau_{0},1_{y}]=[\tau_{1},\xi(\tau_{0},i)(y)]$ by the equivalence relation defining $G(\xi)$. So
\[\begin{array}{ll}
[\tau_{0},1_{y}][\sigma_{0},1_{x}]&=[\tau_{1},\xi(\tau_{0},i)(y)][\sigma_{0},1_{x}]\\
&=[\tau_{1}\sigma_{0},\xi(\tau_{1},\sigma_{0})(x)(\xi(\tau_{0},i)(y)\cdot\sigma_{0})]
\end{array}\]
By (17) in 8.4, we have
\[\xi(\tau_{1},\sigma_{0})=\xi(j \tau_{1},\sigma_{0})=\xi(\tau_{0}i,\sigma_{0}),\]
so using the cocycle condition for $\overset{\tau_{0}}{\gets}\overset{i}{\gets}\overset{\sigma_{0}}{\gets}$,
\[\begin{array}{ll}
\xi(\tau_{1},\sigma_{0})(\xi(\tau_{0},i)\cdot\sigma_{0})&=\xi(\tau_{0}i,\sigma_{0})(\xi(\tau_{0},i)\cdot\sigma_{0})\\
&=\xi(\tau_{0},i\sigma_{0})\xi(i,\sigma_{0})\qquad\mbox{(cocycle cond.)}\\
&=\xi(\tau_{0},i \sigma_{0})\qquad\mbox{(normality)}\\
&=\xi(k\tau_{0}, i\sigma_{0})\qquad\mbox{(by (17))}\\
&=\xi(\tau,\sigma)
\end{array}\]
Thus
\[[\tau_{0},1_{y}][\sigma_{0},1_{x}]=[\tau_{1}\sigma_{0},\xi(\tau,\sigma)(x)].\]
Next,
\[[\tau_{1}\sigma_{0},1_{x}]^{-1}=[(\tau_{1}\sigma_{0})^{-1},\xi(\tau_{1}\sigma_{0},
(\tau_{1}\sigma_{0})^{-1})(x)^{-1}]\]
(cf. 8.5), so we can rewrite definition (23) as
\[\begin{array}{ll}
\zeta(\tau,\sigma)(x)&=[(\tau_{1}\sigma_{0})^{-1},\xi(\tau_{1}\sigma_{0},
(\tau_{1}\sigma_{0})^{-1})(x)^{-1}][\tau_{1}\sigma_{0},\xi(\tau,\sigma)(x)]\\
&=[1_{U},\xi((\tau_{1}\sigma_{0})^{-1},\tau_{1}\sigma_{0})(x)(\xi(\tau_{1}\sigma_{0},
(\tau_{1}\sigma_{0})^{-1})(x)^{-1}\cdot(\tau_{1}\sigma_{0})(x))\xi(\tau,\sigma)(x)]\\
&=[1_{U},\xi(\tau,\sigma)(x)],
\end{array}\]
the latter by (19). This last expression is the image of $\xi(\tau,\sigma)(x)$ under $K\to G(\xi)$, so we have
proved that $\zeta(\tau,\sigma)=\xi(\tau,\sigma)$. This shows that, constructing a groupoid $G(\xi)$, we get
$\zeta=\xi$ back if we choose appropriate liftings in the definition of $\zeta$. It follows that $[\xi]=[\zeta]$,
independent of the liftings.

\subsection{Remark.} Let $E$ be an \'etale groupoid over $T$, as before. Consider the constant circle bundle of
Lie groups $S^{1}\times T\to T$, with \textit{trivial} action of $E$. We denote this bundle simply by $S^{1}$.
Extensions
\[S^{1}\to G\to E\]
of an \'etale groupoid $E$ correspond to bundle gerbes over $E$, as considered in \cite{LupUrb}. Let
$T'\twoheadrightarrow T$ be a surjective \'etale map, and let $E'$ be the \'etale groupoid over $T'$ obtained by
pulling back $E$ along $T'\times T'\to T\times T$. Thus $E'$ is Morita equivalent to $T$. In \cite{KMRW}, Kumjian
et al. consider the colimit $Br(E)$ of extension groups over all $T'\twoheadrightarrow T$, in our notation
\begin{equation}\label{8rem1}
Br(E)=colim_{T'}(\pi_{0}Ext_{T'}(E',S^{1})).
\end{equation}
(Recall that the action on $S^{1}$ is taken to be trivial, and the notation $\mu$ for the action occurring in
Theorem 8.2 is deleted here.) They prove an isomorphism
\begin{equation}\label{8rem2}
Br(E)=H^{2}(E,\underline{S^{1}}).
\end{equation}
This isomorphism is in fact a special case of Theorem 8.2. Indeed, first of all, notice that the sheaf
$\underline{S^{1}}$ of circle valued functions satisfies the conditions of Theorem 7.4 (by the long exact sequence
cohomology  for $0\to \ZZ\to\underline{\RR}\to\underline{S^{1}}\to 0$ and the fact that $\ZZ$ (a constant sheaf)
and $\underline{\RR}$ (a fine sheaf) satisfy these conditions. Thus, by Remark 7.8, we have
\[H^{2}(E',T',S^{1})=H^{2}(E,\underline{S^{1}})=H^{2}(E,\underline{S^{1}}),\]
as soon as $H^{i}(T',\underline{S_{1}})=0$ for $i=1,2$. Thus, for any such $T'$, we have
\begin{equation}\label{8rem3}
\pi_{0}Ext_{T'}(E',S^{1})=H^{2}(E,S^{1}).
\end{equation}
This shows that the l.h.s. of (\ref{8rem3}) is independent of $T'$ provided $H^{i}(T',\underline{S^{1}})=0$ for
$i=1,2$, and that the colimit (\ref{8rem1}) in fact stabilizes. The isomorphism (\ref{8rem2}) becomes a special
case of Theorem 8.2.

A similar remark of course applies not only to $S^{1}$ but to any locally constant bundle of abelian Lie groups.

\newpage
\part{Gerbes and non-abelian cohomology}
In this part, we will discuss the relation between extensions of Lie groupoids and gerbes, with the aim of
describing how the obstruction class of Proposition 5.1 and the cohomological description of extensions in Theorem
8.2 can be adapted to the non-abelian case. This adaption involves non-abelian cohomology in degree 2, for
manifolds as in Proposition 5.1, and for pairs of small categories $(\KE,\KE_{u})$ as in \S7. Both are special
cases of Giraud's general non-abelian cohomology for sites, but each of the two cases here is significantly
simpler. The (first) case is very clearly exposed in \cite{Breen}, and the necessary definitions will be
recapitulated here in a suitable form in \S9. In addition, we state a lemma characterizing neutral gerbes. The
obstruction class is then described in \S10. The (second) case of pairs of small categories is discussed in \S11.

\section{Review of gerbes on manifolds}
We begin by introducing a category of gerbes, somewhat more restrictive than the usual one (cf. Remark 9.3).

\subsection{Sheaves of groupoids.}
Let $M$ be manifold. (In fact, everything in this section applies equally to topological spaces.) a \textit{sheaf
of groupoids} $G$ is a Lie groupoid $G_{1}$ over another manifold $G_{0}$,
\[(\xymatrix{G_{1}\ar@<.5ex>[r]^{s}\ar@<-.5ex>[r]_{t}&G_{0}\ar[r]^{u}&G_{1}},\mbox{ etc.})\]
equipped with \'etale maps (local diffeomorphisms) $p_{0}:G_{0}\to M$ and $p_{1}:G_{1}\to M$, which commute with
all the structure maps of $G$. It follows that $G$ is itself an \'etale groupoid (1.3(f)). The fiber $G_{x}$ of
$G$ over a point $x\in M$ is the discrete groupoid
$\xymatrix{p_{1}^{-1}(x)\ar@<.5ex>[r]\ar@<-.5ex>[r]&p_{0}^{-1}(x)}$. This is a full subgroupoid of $G$. We will
write $G$, or
\[p:G\to M,\]
to denote such a sheaf of groupoids over $M$.

If $p:G\to M$ and $q:H\to M$ are two sheaves of groupoids over $M$, then a map between them is a homomorphism
$\phi:H\to G$ of Lie groupoids such that $p \phi=q$ (i.e. $p_{i}\phi=q_{i}: H_{i}\to M$ for $i=0,1$). We will
continue to use the terminology of \S1. In particular, such a map $\phi$ is called a \textit{weak equivalence} if
it satisfies the conditions of 1.6. In this context, this means that the map $\phi$ induces for each point $x\in
M$ an ordinary equivalence of small categories $H_{x}\to G_{x}$ between the fibers. Two sheaves of groupoids
$G_{1}$ and $G_{2}$ over $M$ are said to be (Morita) equivalent if there are two such weak equivalences
$G_{1}\gets H\to G_{2}$ over $M$.

\subsection{Gerbes.}
A sheaf of groupoids $G$ over $M$ is called a \textit{gerbe} if each stalk $G_{x}$ is a non-empty transitive
(discrete) groupoid. In other words, $G$ is a gerbe iff the \'etale maps $p_{0}:G_{0}\to M$ and $(s,t):G_{1}\to
G_{0}\times_{M}G_{0}$ are both surjective. This condition is invariant under Morita equivalence.

\subsection{Remark.}
This Definition 9.2 is more restrictive than the usual one (\cite{Giraud, DeMi, Breen}), but it suffices for our
purposes, essentially because we only need to consider gerbes whose band (cf. 9.5 below) is isomorphic to a sheaf
of groups.

\subsection{Outer isomorphisms.}
Suppose $\underline{K}$ and $\underline{L}$ are sheaves of groups over $M$. We write
$Iso(\underline{K},\underline{L})$ for the sheaf of isomorphisms from $\underline{K}$ to $\underline{L}$. Its
stalk at $x\in M$ consists of germs of isomorphisms $\underline{K}|_{U}\to\underline{L}|_{U}$ where $U$ ranges
over nbds of $x$. We write $Out(\underline{K},\underline{L})$ for the quotient sheaf of outer isomorphisms. Thus,
a section of $Out(\underline{K},\underline{L})$ over $U$ is \textit{represented} by a cover $U=\bigcup U_{i}$, and
isomorphisms $\phi_{i}:\underline{K}|_{U_{i}}\to\underline{L}|_{U_{i}}$, such that for any two indices $i$ and
$j$, there is a cover $U_{i j}=\bigcup V_{i j}^{\xi}$ for which $\phi_{i}|_{V_{i j}^{\xi}}$ and $\phi_{j}|_{V_{i
j}^{\xi}}$ differ by conjugation by a section $\lambda_{i j}^{\xi}$ of $\underline{L}$ over $V_{i j}^{\xi}$. Two
such families $\{\phi_{i}:\underline{K}|_{U_{i}}\to \underline{L}|_{U_{i}}\}$ and
$\{\psi_{j}:\underline{K}|_{W_{j}}\to\underline{L}|_{W_{j}}\}$ represent the \textit{same} section of
$Out(\underline{K},\underline{L})$ if there is a common refinement, say $O_{i j}\subset U_{i}\cap W_{j}$, over
which they differ  by conjugation by a section of $\underline{L}$. A section of $Out(\underline{K},\underline{L})$
over $U$ will be referred to as an ``\textit{outer isomorphism}'' $\underline{K}|_{U}\to\underline{L}|_{U}$. Note
that if $\underline{K}$ and $\underline{L}$ are abelian, then
$Out(\underline{K},\underline{L})=Iso(\underline{K},\underline{L})$.

\subsection{Bands.}
(cf. \cite{Giraud, DeMi}) A band over $M$ is represented by a cover $M=\bigcup U_{i}$ and for each $i$ a sheaf of
groups $\underline{K}_{i}$ on $U_{i}$, together with outer isomorphisms
\[\lambda_{i j}:\underline{K}_{j}|_{U_{i j}}\overset{\sim}{\to}\underline{K}_{i}|_{U_{i j}}\]
satisfying the cocycle condition $\lambda_{i j}\lambda_{j k}=\lambda_{i k}$ (equality as outer isomorphisms!). If
$\{V_{i}\}$ is a refinement of $\{U_{i}\}$, say by $V_{j}\subset U_{\alpha(j)}$, then the restricted data,
$\underline{K}_{j}=\underline{K}_{\alpha(j)}|_{V_{j}}$ and $\lambda_{j j'}=\lambda_{\alpha(j)\alpha(j')}|_{V_{j
j'}}$, are taken to represent the \textit{same} band.

If $\underline{K}$ and $\underline{L}$ are two bands over $M$, we may assume that they are represented over the
same cover $\{U_{i}\}$. An isomorphism $\underline{K}\to \underline{L}$ is then defined to be a system of outer
isomorphisms $\theta_{i}:\underline{K}_{i}\to \underline{L}_{i}$, which are compatible with the cocycles
$\lambda_{i j}$ for $\underline{K}$ and $\mu_{i j}$ for $\underline{L}$, in the sense that $\theta_{i}\lambda_{i
j}=\mu_{i j}\theta_{j}$ as outer isomorphisms over $U_{i j}$. Passing to finer cover defines the \textit{same}
morphism.

If $\underline{K}=(\underline{K}_{i},\lambda_{i j})$ is a band and each $\underline{K}_{i}$ is abelian, then the
$\lambda_{i j}$ satisfy the cocycle condition on the nose, and there is a sheaf of abelian groups $\underline{K}$
on $M$, such that $\underline{K}_{i}\cong \underline{K}|_{U_{i}}$ (by an isomorphism compatible with the
$\lambda_{i j}$). A similar remark applies to (iso)morphisms. Thus, the category of abelian bands is equivalent to
that of abelian sheaves.

\subsection{The band of a gerbe.} (1) Let $p:G\to M$ be a gerbe over $M$. Since $p_{0}:G_{0}\to M$ is surjective,
there exists an open cover $M=\bigcup U_{i}$ over which there are sections
\[s_{i}:U_{i}\to G_{0}.\]
Let $\underline{Aut}(s_{i})$ be the automorphism group --- this is a sheaf of groups on $U_{i}$ with as stalk at
$x\in U_{i}$ the isotropy group of $s_{i}(x)$ in the groupoid $G_{x}$. Since $(s,t):G_{1}\to G_{0}\times_{M}G_{0}$
is surjective, there is a cover $U_{i j}=\bigcup V_{i j}^{\xi}$ over which there are sections $\beta_{i
j}^{\xi}:V_{i j}^{\xi}\to G_{1}$ such that
\[\beta_{i j}(x):s_{j}(x)\to s_{i}(x)\]
in $G_{x}$. Conjugation by $\beta_{i j}^{\xi}$ defines an isomorphism $(\beta_{i j}^{\xi})_{*}:
\underline{Aut}(s_{j})\to \underline{Aut}(s_{i})$ over $V_{i j}^{\xi}$. Since any two $\beta_{i j}^{\xi}$ and
$\beta_{i j}^{\eta}$ differ by conjugation by a section of $\underline{Aut}(s_{i})$ over $V_{i j}^{\xi \eta}=V_{i
j}^{\xi}\cap V_{i j}^{\eta}$, these $(\beta_{i j}^{\xi})_{*}$ patch together to an \textit{outer} isomorphism
\[\lambda_{i j}:\underline{Aut}(s_{j})|_{U_{i j}}\to \underline{Aut}(s_{i})|_{U_{i j}}.\]
Thus we have defined a band, called the band of the gerbe $G$ and denoted
\[Band(G).\]

(2) It often happens that one can find sections $\beta_{i j}: U_{i j}\to G_{1}$ with $\beta_{i j}(x):s_{j}(x)\to
s_{i}(x)$, so that it isn't necessary to pass to a cover $V_{i j}^{\xi}$ of $U_{i j}$. For example, this is always
the case if $H^{1}(U_{i j},\underline{Aut}(s_{i}))=0$, by an argument similar to the one in the proof of
Proposition 5.1 above. This can often be arranged by choosing the $\{U_{i}\}$ to be a ``good'' cover, so that the
$U_{i j}$ are contractible. We will often \textit{assume} that such $\beta_{i j}$ can be found, just to simplify
the notation and avoid too many indices. However this assumption is never essential.

(3) Let us replace the sections $s_{i}$ in (1) above by sections $r_{i}$, and let us suppose (cf.(2)) that there
are sections $\alpha_{i}:U_{i}\to G_{1}$ with $\alpha_{i}(x): r_{i}(x)\to s_{i}(x)$ in $G_{x}$. The conjugation by
$\alpha_{i}$ defines a map
\[(\alpha_{i})_{*}:\underline{Aut}(r_{i})\to\underline{Aut}(s_{i})\]
and these patch together to a \textit{canonical} isomorphism of bands $\alpha :\{\underline{Aut}(r_{i})\}
\to\{\underline{Aut}(s_{i})\}$. Strictly speaking, one should view $Band(G)$ as the system of these bands and
canonical isomorphisms between them.

(4) If $\underline{L}$ is an abstractly given band, a \textit{gerbe with band} $\underline{L}$ is a gerbe $G$
together with an isomorphism of bands $\theta:Band(G)\to \underline{L}$. If $Band(G)=\{Aut(s_{i}),\lambda_{i j}\}$
as in (1), and the cover $\{U_{i}\}$ is sufficiently fine, then $\underline{L}$ can be represented by
$(\underline{L}_{i},\mu_{i j})$ for this cover, and $\theta$ by
$\theta_{i}:\underline{Aut}(s_{i})\to\underline{L}_{i}$. If $\alpha:\{\underline{Aut}(r_{i})\}\to
\{\underline{Aut}(s_{i})\}$ is a canonical isomorphism as above, then $\theta:\{\underline{Aut}(s_{i})\} \to
\underline{L}$ should be identified with $\theta\alpha:\{\underline{Aut}(r_{i})\}\to \underline{L}$.

\subsection{Non-abelian $H^{2}$.}
(\cite{Giraud}) Let $\underline{K}$ be a band on $M$, and let $(G\overset{p}{\to}M,
\theta:Band(G)\to\underline{K})$ and $(H\overset{q}{\to}M, \tau:Band(H)\to\underline{K})$ be two gerbes with band
$\underline{K}$. We consider only those weak equivalences $\phi:H\to G$ over $M$ for which the induced map
$Band(\phi):Band(H)\to Band(G)$ satisfies $\theta\circ Band(\phi)=\tau$ (as outer morphisms between bands, cf.
9.5). The same applies to Morita equivalences $G_{1}\gets H\to G_{2}$. The set of equivalence classes of gerbes
with band $\underline{K}$ is denoted
\[H^{2}(M,\underline{K})\]
It agrees with the ordinary sheaf cohomology if $\underline{K}$ is abelian.

If $\underline{K}$ is a sheaf of groups, it can be viewed as a band (represented on the trivial cover consisting
of just $M$ itself), and also as a gerbe (namely the gerbe $\xymatrix{G_{1}\ar@<.5ex>[r]\ar@<-.5ex>[r]&G_{0}}$
with $G_{1}=\underline{K}$ and $G_{0}=M$). The identity $\underline{K}\to \underline{K}$ gives $\underline{K}$ the
structure of a gerbe with band $\underline{K}$. Thus $H^{2}(M,\underline{K})$ is a \textit{pointed set}, with base
point the equivalence class of $(\underline{K},id)$. A gerbe with band $\underline{K}$ is called \textit{neutral}
if it is equivalent to $(\underline{K},id)$.

We now characterize neutral gerbes, in Lemma 9.9 below. To simplify the notation, we will only do this under a
hypothesis as in 9.6(2). However, a similar characterization is valid without this hypothesis. First we need a
preliminary remark.

\subsection{Remark.}
Let $(G,\theta)$ be a gerbe with band a sheaf of groups $\underline{K}$. Let $\theta:Band(G)\to \underline{K}$ be
represented by isomorphisms
\[\theta_{i}:\underline{Aut}(s_{i})\to \underline{K}|_{U_{i}}.\]
Under the simplifying hypothesis 9.6(2), there are sections $\beta_{i j}:s_{j}\to s_{i}$ over $U_{i j}$, for which
the diagram
\begin{equation}\label{rr10}
\xymatrix{\underline{Aut}(s_{j})|_{U_{i j}}\ar[r]^-{\theta_{j}}\ar[d]_-{(\beta_{i j})_{*}}& \underline{K}|_{U_{i
j}}\ar@{=}[d]\\ \underline{Aut}(s_{i})|_{U_{i j}}\ar[r]^-{\theta_{i}}& \underline{K}|_{U_{i j}}}
\end{equation}
commutes as outer isomorphisms. Thus, it commutes locally up to conjugation by local sections $g_{i j}^{\xi}$ of
$\underline{Aut}(s_{i})$ on a cover $V_{i j}^{\xi}$ of $U_{i j}$. Under similar simplifying hypotheses (that
$H^{1}(U_{i j}, Z(\underline{Aut}(s_{i}))=0$), we may assume that there is in fact just \textit{one} section $g_{i
j}$ of $\underline{Aut}(s_{i})$ such that $\theta_{i}(g_{i j})_{*}(\beta_{i j})_{*}=\theta_{j}$. Thus, first
choosing \textit{actual} maps $\theta_{i}$ and $\theta_{j}$ to represent $\theta$, and then replacing the arrow
$\beta_{i j}(x): s_{j}(x)\to s_{i}(x)$ in $G_{x}$ by $g_{i j}(x)\beta_{i j}(x)$, we may assume that the square (1)
commutes in the category of sheaves of groups, and not only as outer isomorphisms. Given such actual maps
$\theta_{i}$, $\theta_{j}$ and $(\beta_{i j})_{*}$ making (1) commute, \cite{Breen} represents the gerbe $G$ by
the ``non-abelian 2-cocycle'' $\{\lambda_{i j}, g_{i j k}\}$, where $\lambda_{i j}=\phi_{i}(\beta_{i
j})_{*}\phi_{j}^{-1}$, and $g_{i j k}=\theta_{i}(\beta_{i j} \beta_{j k}\beta_{i k}^{-1})\in \Gamma(U_{i j
k},\underline{K})$. For a \textit{sheaf} of groups $\underline{K}$, we see that we can take $\lambda_{i j}=1$, and
the gerbe $G$ is represented just by the cocycle $g_{i j k}$.

\subsection{Lemma.}
\textit{Let $\underline{K}$ be sheaf of groups on $M$, and let $(G,\theta:
Band(G)\overset{\sim}{\to}\underline{K})$ be a gerbe with band $\underline{K}$. Then $G$ is neutral iff there are
sections $s_{i}$ and $\beta_{i j}:s_{j}\to s_{i}$ as above, such that $\theta$ can be represented by an actually
commutative diagram (27) and the cocycle condition $\beta_{i j}\beta_{j k}=\beta_{i k}$ is satisfied.}

\begin{proof}
($\Rightarrow$) Suppose we are given equivalences
\[G\overset{\phi}{\gets}H\overset{\psi}{\to}\underline{K}\]
compatible with $\theta$. Choose local sections $r_{i}: U_{i}\to H_{0}$ of $H_{0}\to M$ on a cover $\{U_{i}\}$ of
$M$. Since $\psi$ is an equivalence, there is a unique section $\alpha_{i j}$ of $H_{1}$ over $U_{i j}$ such that
$\alpha_{i j}(x):r_{j}(x)\to r_{i}(x)$ in $H_{x}$ is mapped to the identity in $\underline{K}_{x}$. Then the
$\alpha_{i j}$ satisfy the cocycle condition. Let $s_{i}=\phi(r_{i})$ and $\beta_{i j}=\phi(\alpha_{i j})$, and
represent $\theta_{i}$ by
\[\underline{Aut}(s_{i})\overset{\phi^{-1}}{\to}\underline{Aut}(r_{i})\overset{\psi}{\to}\underline{K}.\]

($\Leftarrow$) Suppose we are given $s_{i}$, $\beta_{i j}$ and actual maps $\theta_{i}: \underline{Aut}(s_{i})\to
\underline{K}|_{U_{i}}$ as in the statement of the lemma. To define a sheaf of groupoids $q: H\to M$, let
$H_{0}=\coprod U_{i}$ with $q_{0}:H_{0}\to M$ the evident map. Let $\phi_{0}: H_{0}\to G_{0}$ be the coproduct of
the maps $s_{i}: U_{i}\to G_{0}$, and let $H_{1}$ and $\phi_{1}$ be defined by the pullback
\[\xymatrix{H_{1}\ar[r]^{\phi_{1}}\ar[d]_{(s,t)}&G_{1}\ar[d]^{(s,t)}\\H_{0}\times H_{0}\ar[r]^{\phi_{0}\times\phi_{0}}&
G_{0}\times G_{0}}.\] We can write a point of $H_{0}$ as a pair $(x,i)$ where $x\in U_{i}$. Arrows $g:(x,i)\to
(y,j)$ in $H$ are arrows $s_{i}(x)\to s_{j}(y)$ in $G$; these only exist if $x=y$. Define $\psi: H\to
\underline{K}$, for an arrow $g:(x,i)\to (x,j)$, by
\[\psi(g)=\theta_{i}(\beta_{i j} g).\]
Then $\psi$ is a homomorphism. Indeed, using that (27) commutes and that $\beta_{i j}\beta_{j k}=\beta_{i k}$, we
find for $h:(x,j)\to (x,k)$ that
\[\begin{array}{ll}
\psi(h)\psi(g)&=\theta_{j}(\beta_{j k} h)\theta_{i}(\beta_{i j}g)\\
&=\theta_{i}(\beta_{i j}\beta_{j k} h\beta_{i j}^{-1})\theta_{i}(\beta_{i j}g)\\
&=\theta_{i}(\beta_{i k}h g)\\
&=\psi(h g).
\end{array}\]
The resulting diagram of equivalences $G\overset{\phi}{\gets}H\overset{\psi}{\to}\underline{K}$ is compatible with
$\theta$ and shows that $G$ is neutral.
\end{proof}

\section{The non-abelian obstruction class}
Let us go back to the context of Section 5, and consider the restriction functor
\begin{equation}\label{rr11}R:Ext_{M}(E,K)\to Ext_{T}(E_{T},K_{T}).\end{equation}
Here $E$ is a foliation groupoid over $M$, and $K$ is a bundle of Lie groups over $M$. Furthermore, $T$ is a
complete transversal to the orbit foliation of $E$, while $E_{T}$ and $K_{T}$ are the restriction of $E$ and $K$
to $T$. Recall that $E_{T}$ is an \'etale groupoid.

Consider an extension of Lie groupoids over $T$,
\begin{equation}\label{rr12}K_{T}\to B\to E_{T}.\end{equation}
We wish to formulate an obstruction to the problem of finding an extension $K\to G\to E$ over $M$ which restricts
to an extension isomorphic to $B$ over $T$,
\[\xymatrix{K_{T}\ar@{=}[d]\ar[r]&G_{T}\ar[r]\ar[d]^{\sim}& E_{T}\ar@{=}[d]\\
K_{T}\ar[r]&B\ar[r]&E_{T}}\] In \S5, we saw that this obstruction lies in $H^{2}(M,\underline{K})$ in case
$\underline{K}$ is abelian, and we made explicit use of the action of $E$ on $K$. In the present non-abelian case,
we will first formulate the ``action'' of $E$ on $K$ as an isomorphism of bands.

\subsection{The band $\underline{K}^{\gamma}$.}
Let $M=\bigcup U_{i}$ be an open cover for which there exist submersions $\gamma_{i}:U_{i}\to T$ and smooth maps
$\tau_{i}:U_{i}\to E$ with $\tau_{i}(x):x\to\gamma_{i}(x)$ for $x\in U_{i}$, as in the proof of Proposition
4.5(i). Let us write $\tau_{i j}=\tau_{i}\tau_{j}^{-1}:\gamma_{j}(x)\to \gamma_{i}(x)$. This is a cocycle with
values in $E_{T}$, by help of which one can reconstruct $E$ from $E_{T}$; see Remark 4.6. Let
$\underline{K}_{i}^{\gamma}$ be the sheaf of sections of the pullback bundle $\gamma_{i}^{*}(K_{T})$. So
$\underline{K}_{i}^{\gamma}$ is a sheaf of groups on $U_{i}$. Suppose we can lift the arrows $\tau_{i j}(x)$ in
$E_{T}$ to $\beta_{i j}(x):\gamma_{j}(x)\to \gamma_{i}(x)$ in $B$. Then conjugation by $\beta_{i j}$ defines an
isomorphism
\[(\beta_{i j})_{*}:\underline{K}_{j}^{\gamma}\to\underline{K}_{i}^{\gamma}\]
of sheaves of groups. Thus, these $\underline{K}_{i}^{\gamma}$ and $(\beta_{i j})_{*}$ together define a band
$\underline{K}^{\gamma}$. The same is true if we can only lift the $\tau_{i j}$ locally to $\beta_{i j}^{\xi}$,
since the $(\beta_{i j})_{*}$ need only to be given as outer isomorphisms.

\subsection{The groupoid $Ext^{\theta}(E,K)$.}
Suppose that the extension $K_{T}\to B\to E_{T}$ is the restriction of an extension $K\to G\to E$. Suppose the
arrows $\tau_{i}(x):x\to \gamma_{i}(x)$ can be lifted to $\beta_{i}(x):x\to \gamma_{i}(x)$. Then conjugation with
$\beta_{i}(x)$ defines an isomorphism $\theta_{i}:\underline{K}_{i}^{\gamma}\to\underline{K}$, not depending on
$\beta_{i}$ when viewed as an outer isomorphism. These $\theta_{i}$ patch together to an isomorphism of bands
\[\theta:\underline{K}^{\gamma}\to\underline{K}.\]
This isomorphism is what replaces the action of $E$ on $K$ in the abelian case. It should be compatible with the
already existing similar isomorphism over $T$.

In detail, let $\underline{K}_{T}$ be the sheaf of sections of $K_{T}$. This sheaf is not isomorphic to the
restriction of the sheaf $\underline{K}$ to $T$, but there is an obvious restriction map $R:\underline{K}|_{T}\to
\underline{K}_{T}$. Similarly, let $\underline{K}_{T,i}^{\gamma}$ be the sheaf of sections on $U_{i}\cap T$ of the
pullback of $K_{T}$ along $\gamma_{i}:U_{i}\cap T\to T$. These sheaves $\underline{K}_{T,i}^{\gamma}$ form a band
$\underline{K}^{\gamma}_{T}$ on $T$, and there is a restriction map $R^{\gamma}:\underline{K}^{\gamma}|_{T}\to
\underline{K}^{\gamma}_{T}$. The extension $K_{T}\to B\to E_{T}$ induces an isomorphism
$\theta_{i}^{T}:\underline{K}_{T,i}^{\gamma}\to\underline{K}_{T}$ over $U_{i}\cap T$, and these patch together to
an isomorphism of bands $\theta^{T}:\underline{K}_{T}^{\gamma}\to\underline{K}_{T}$. The result is a diagram
\begin{equation}\label{rr13}
\xymatrix{\underline{K}^{\gamma}|_{T}\ar[r]^-{\theta|_{T}}_{\sim}\ar[d]_{R^{\gamma}}&
\underline{K}|_{T}\ar[d]^{R}\\\underline{K}_{T}^{\gamma}\ar[r]^{\theta^{T}}_{\sim}&\underline{K}_{T}}
\end{equation}
which commutes up to inner automorphisms.

Given such an isomorphism of bands $\theta:\underline{K}^{\gamma}\to\underline{K}$, we write
$Ext^{\theta}_{M}(E,K)$ for the full subgroupoid of $Ext_{M}(E,K)$, consisting of extensions which induce the
isomorphism $\theta:\underline{K}^{\gamma}\to\underline{K}$. (This groupoid plays the same r\^ole as
$Ext^{\mu}_{M}(E,K)$ in the abelian case; cf. 6.8. The functor $R$ in (28) now restricts to a functor
\begin{equation}\label{rr14}
R: Ext_{M}^{\theta}(E,K)\to Ext_{T}^{\theta^{T}}(E_{T},K_{T}).
\end{equation}
Given an extension $K_{T}\to B\to E_{T}$ inducing $\theta^{T}$, and a ``lifting''
$\theta:\underline{K}^{\gamma}\to\underline{K}$ making (30) commute, we will now define an obstruction for
$K_{T}\to B\to E_{T}$ to be in the image of $R$ in (28), cf. Proposition 10.4 below.

\subsection{The gerbe $\hat{B}$.}
From an extension $K_{T}\to B\to E_{T}$ as in (29) we construct a gerbe $\hat{B}$ on $T$, as follows. Let
\[\hat{B}_{0}=E\times_{M}T=\{e:x\to y\mbox{ in }E\quad|\quad y\in T\},\]
and let $p:\hat{B}_{0}\to M$ be the map sending $e:x\to y$ to $x$. Then $p$ is an \'etale map because $T$ is a
transversal to the orbits of $E$. Consider also the space $\hat{B}_{0}\times_{T}B$ of pairs $(e:x\to y, b:y \to
z)$ where $e\in E$, $b\in B$ and $y,z\in T$. There is an evident surjective submersion
\[\mu:\hat{B}_{0}\times_{T}B\to \hat{B}_{0}\times_{M}\hat{B}_{0},\]
mapping $(e:x\to y,b:y\to z)$ to $(e,\pi(b)e)$. Let $\hat{B}_{1}$ be the sheaf of sections of $\mu$. Then
$\hat{B}_{1}$ is a sheaf over $\hat{B}_{0}\times_{M}\hat{B}_{0}$, the \'etale projection of which we denote by
\[(s,t):\hat{B}_{1}\to \hat{B}_{0}\times_{M}\hat{B}_{0}.\]
Define $p_{1}:\hat{B}_{1}\to M$ as the composition of $(s,t)$ with $\hat{B}_{0}\times_{M}\hat{B}_{0}\to M$. Then
$\xymatrix{\hat{B}_{1}\ar@<.5ex>[r]^{s}\ar@<-.5ex>[r]_{t}&\hat{B}_{0}}$ has the structure of a sheaf of groupoids
over $M$.

The stalk at $x\in M$ is the discrete groupoid whose objects are such $e:x\to y$ in $E$ with $y\in T$. Since
$\hat{B}_{0}\to M$ is \'etale, there is unique germ of a local section $\tilde{e}:U_{x}\to \hat{B}_{0}$ with
$\tilde{e}(x)=e$; write $\tilde{e}(x'):x'\to y'$, where $x'\in U_{x}$. Similarly, let $\tilde{f}$ be a local
section through $f:x\to z$. Then an arrow from $(e:x\to y)$ to $(f:x\to z)$ in $\hat{B}$ is the germ of a map
$U_{x}\overset{\tilde{b}}{\to}B$, such that $\tilde{b}(x'):y'\to z'$ is an arrow in $B$ with
$\pi(\tilde{b}(x'))\tilde{e}(x')=\tilde{f}(x')$ in the groupoid $E$.

It is clear that $\hat{B}$ thus defined is a gerbe. Indeed, $\hat{B}_{0}\to M$ is surjective because the
transversal $T$ is complete, and $\hat{B}_{1}\overset{(s,t)}{\to}\hat{B}_{0}\times_{M}\hat{B}_{0}$ is surjective
because it is the sheaf of sections of the surjective submersion $\hat{B}_{0}\times_{T}B
\to\hat{B}_{0}\times_{B}\hat{B}_{0}$.

The $\tau_{i}(x):x\to \gamma_{i}(x)$ for $x\in U_{i}$ (cf. 10.1) define sections $\tau_{i}: U_{i}\to\hat{B}_{0}$,
and $\underline{Aut}(\tau_{i})=\underline{K}_{i}^{\gamma}$. So $\hat{B}$ is a gerbe with band
$\underline{K}^{\gamma}$; in fact, if we choose these sections $\tau_{i}$ to represent $Band(\hat{B})$, then the
required isomorphism $Band(\hat{B})\to\underline{K}^{\gamma}$ is the identity. By composition with
$\theta:\underline{K}^{\gamma}\to \underline{K}$, we thus find a gerbe $(\hat{B}, \theta:Band(\hat{B})\to
\underline{K})$ representing an element of $H^{2}(M,\underline{K})$. We claim that this is the obstruction class:

\subsection{Proposition.}
\textit{The extension $K_{T}\to B\to E_{T}$ is in the image of the restriction functor (28) (up to isomorphism),
iff the class of $(\hat{B}, \theta:Band(\hat{B})\to \underline{K})$ vanishes in $H^{2}(M,\underline{K})$.}

(This is the non-abelian analogue of Proposition 5.1. It implies Corollaries similar to Corollary 5.1 and Theorem
6.6.)
\begin{proof}
(Sketch) By Lemma 9.9, $\hat{B}$ is neutral iff there is a cocycle $\beta_{i j}:\tau_{j}\to \tau_{i}$ with values
in $B$, i.e. $\beta_{i j}(x):\tau_{j}(x)\to \tau_{i}(x)$ for $x\in U_{i}$ with $\beta_{i j}(x)\beta_{j
k}(x)=\beta_{i k}(x)$. This is equivalent to finding a cocycle $\beta_{i j}$ which lifts $\tau_{i j}$, and this is
precisely the obstruction to enlarging the extension $K_{T}\to B\to E_{T}$ over $T$ to an extension $K\to G\to E$
over $M$; cf. the proof of Proposition 5.1 and Remark 4.6.
\end{proof}

\section{Relative non-abelian cohomology of small categories}
The purpose of this section is to introduce a relative non-abelian cohomology
$H^{2}(\KE,\mathcal{D},\underline{K})$ for a pair of small categories $\mathcal{D}\subset\KE$ and suitable
coefficients in $\underline{K}$. This will enable us to formulate an analogue of Theorem 8.2 in the non-abelian
context.

We begin by considering the non-relative case. Throughout this section, we will use the notation
\[g_{*}:G\to G\]
for the inner automorphism $g_{*}(x)=g x g^{-1}$, for any element $g$ in a group $G$.

\subsection{Presheaves and bands.}
Let $\KE$ be a small category. A \textit{presheaf of groups} on $\KE$ is a functor $\underline{K}: \KE^{op}\to
Groups$. Thus, $\underline{K}$ assigns to each object $U$ of $\KE$ a group $\underline{K}(U)$, and to each arrow
$\sigma:V\to U$ a group homomorphism which we denote as $\underline{K}(\sigma)$ or $\sigma^{*}:
\underline{K}(U)\to \underline{K}(V)$, or as the map $k\mapsto k\cdot\sigma$ for $k\in \underline{K}(U)$ (as in
\S7). A map between two such presheaves is by definition a natural transformation. In a similar way, one defines
presheaves on $\KE$ with values in an arbitrary category. A \textit{band} on $\KE$ is a presheaf $\underline{K}$
on $\KE$ with values in the category of groups and outer homomorphisms (Thus, for arrow
$W\overset{\tau}{\to}V\overset{\sigma}{\to}U$, the functoriality condition $(\sigma\tau)^{*}=\tau^{*}\sigma^{*}$
holds only up to conjugation by an element of $\underline{K}(W)$). Any presheaf of groups can be viewed as a band.

\subsection{Presheaves of groupoids and bands.}
Let $\mathcal{G}$ be a presheaf on $\KE$ with values in non-empty transitive discrete groupoids. Choose for each
object $U$ in $\KE$ an object $x_{U}$ of $\mathcal{G}(U)$. Next, choose for each arrow $\sigma:V\to U$ in $\KE$ an
arrow
\begin{equation}\label{rr15}\tilde{\sigma}:x_{V}\to\sigma^{*}(x_{U})\end{equation}
in $\mathcal{G}(V)$ (there is at least one such arrow because $\mathcal{G}(V)$ is transitive). Then let
$\underline{Aut}(x_{U})=\mathcal{G}(U)(x_{U},x_{U})$ be the isotropy group of $x_{U}$ in the groupoid
$\mathcal{G}(U)$, and define
\[\tilde{\sigma}^{*}:\underline{Aut}(x_{U})\to \underline{Aut}(x_{V})\]
to be the composition of $\sigma^{*}$ and conjugation with $\tilde{\sigma}$,
\[\tilde{\sigma}^{*}(g)=\tilde{\sigma}^{-1}\sigma^{*}(g)\tilde{\sigma}\]
When viewed as an outer homomorphism, the map $\tilde{\sigma}^{*}$ is independent of the choice of
$\tilde{\sigma}$. Thus we have defined a band on $\KE$, called ``the'' band $\mathcal{G}$ and denoted
$Band(\mathcal{G})$.

Let $\underline{K}$ be \textit{any} band on $\KE$. For the purposes of this section, we define a gerbe on $\KE$
with band $\underline{K}$ as a pair $(\mathcal{G},\theta)$ where $\mathcal{G}$ is a presheaf of non-empty
transitive groupoids and $\theta: Band(\mathcal{G})\to \underline{K}$ is an isomorphism of bands. Here
$Band(\mathcal{G})=\{\underline{Aut}(x_{U}),\tilde{\sigma}^{*}\}$ as above, and depends on the choice of the
$x_{U}$. If $y_{U}\in\mathcal{G}(U)$ is another choice, arrows $\alpha_{U}:y_{U}\to x_{U}$ induce by conjugation a
\textit{canonical} isomorphism of bands $\alpha_{*}:\{\underline{Aut}(y_{U}\}\to\{\underline{Aut}(x_{U)}\}$, and
we identify $(\mathcal{G},\theta)$ with $(\mathcal{G},\theta\alpha_{*})$ where
$\theta_{U}:\underline{Aut}(x_{U})\to \underline{K}(U)$, exactly as in 9.6 above.

More generally \cite{Giraud}, a gerbe on $\KE$ is usually defined as a fibered category over $\KE$ with connected
transitive non-empty groupoids as fibers. But presheaves suffice for our purposes.

\subsection{Cocycles of degree 2.}
Let $\underline{K}$ be a band on $\KE$. A 2-cocycle on $\KE$ with values in $\underline{K}$ is a pair
\[(R,\xi)\]
where $R$ assigns to each $\sigma:V\to U$ in $\KE$ an actual group homomorphism
$R(\sigma):\underline{K}(U)\to\underline{K}(V)$ which represents the conjugacy class $\sigma^{*}:
\underline{K}(U)\to\underline{K}(V)$, and $\xi$ assigns to each composable pair
$U\overset{\sigma}{\gets}V\overset{\tau}{\gets}W$ of arrows in $\KE$ an element
$\xi(\sigma\tau)\in\underline{K}(W)$. Furthermore, there are conditions that should be satisfied: first, the
failure of functoriality of $R$ is measured by conjugation with $\xi$:
\begin{equation}\label{rr16}
\xi(\sigma,\tau)_{*}R(\tau)R(\sigma)=R(\sigma \tau);
\end{equation}
and secondly, $\xi$ satisfies the ``cocycle condition'' for any three composable arrows
$N\overset{\rho}{\gets}U\overset{\sigma}{\gets}V\overset{\tau}{\gets}W$ in $\KE$,
\begin{equation}\label{rr17}
\xi(\rho,\sigma\tau)\xi(\sigma,\tau)=\xi(\rho \sigma,\tau)R(\tau)(\xi(\rho,\sigma))
\end{equation}
(equality in $\underline{K}(W)$).

\subsection{Example.}
Suppose $(\mathcal{G},\theta)$ is a gerbe with band $\underline{K}$, as in 11.2. As there, choose objects $x_{U}$
and arrows $\tilde{\sigma}:x_{V}\to \sigma^{*}(x_{U})$ to represent $Band(\mathcal{G})$; and next choose actual
homomorphisms $\theta_{U}:\underline{Aut}(x_{U})\to \underline{K}(U)$ to represent the outer homomorphism
$\theta_{U}$ (this introduces some ambiguity in notation, but shouldn't lead to confusion). Now choose an actual
homomorphism $R(\sigma):\underline{K}(U)\to\underline{K}(V)$ for each $\sigma:V\to U$ in $\KE$, to represent
$\sigma^{*}:\underline{K}(U)\to\underline{K}(V)$, and with the additional property that the square
\[\xymatrix{\underline{Aut}(x_{U})\ar[r]^-{\theta_{U}}\ar[d]^{\tilde{\sigma}^{*}}&\underline{K}(U)\ar[d]^{R(\sigma)}\\
\underline{Aut}(x_{V})\ar[r]^-{\theta_{V}}&\underline{K}(V)}\] actually commutes (and not just as outer
homomorphisms). For each composable pair $U\overset{\sigma}{\gets}V\overset{\tau}{\gets}W$, we have chosen arrows
$\widetilde{\sigma\tau}:x_{W}\to(\sigma\tau)^{*}(x_{U}) =\tau^{*}\sigma^{*}(x_{U})$, together with
$\tilde{\sigma}:x_{V}\to \sigma^{*}(x_{U})$ and $\tilde{\tau}:x_{W}\to \tau^{*}(x_{V})$, and we write
\[\xi_{0}(\sigma,\tau)=(\widetilde{\sigma\tau)})^{-1}\tau^{*}(\tilde{\sigma})\tilde{\tau}\in \underline{Aut}(x_{W}),\]
and
\[\xi(\sigma,\tau)=\theta_{W}(\xi_{0}(\sigma,\tau)).\]
Then $(R,\xi)$ is a cocycle with values in $\underline{K}$. We omit the straightforward verification of the
identities (33) and (34).

\subsection{Equivalences of cocycles.}
We introduce the following equivalence relation on 2-cocycles on $\KE$ with values $\underline{K}$. One such
cocycle $(R,\xi)$ is equivalent to another one $(S,\zeta)$ if there are \textit{inner} automorphisms
$\mu_{U}:\underline{K}(U)\to \underline{K}(U)$ (for each object $U$ in $\KE$) and elements $f(\sigma)\in
\underline{K}(V)$ (for each arrow $\sigma:V\to U$ in $\KE$) such that
\begin{equation}\label{rr18}
S(\sigma)=f(\sigma)_{*}\mu_{V}R(\sigma)\mu_{U}^{-1}
\end{equation}
\begin{equation}\label{rr19}
\mu_{W}\xi(\sigma,\tau)=f(\sigma\tau)^{-1}\zeta(\sigma,\tau)S(\tau)(f(\sigma))f(\tau)
\end{equation}
These relations are analogous to \cite{Breen}, page 42. It is easier to understand this equivalence relation in
two steps: first one may replace each $R(\sigma)$ by $\mu_{V}R(\sigma)\mu_{U}^{-1}$ and each $\xi(\sigma,\tau)$ by
$\mu_{W}\xi(\sigma,\tau)$, to obtain an equivalent cocycle. Secondly, one may replace $R(\sigma)$ by
$f(\sigma)_{*}R(\sigma)$ and $\xi(\sigma,\tau)$ by
$f(\sigma\tau)\xi(\sigma,\tau)R(\tau)(f(\sigma)^{-1})f(\tau)^{-1}$, to obtain an equivalent cocycle. The set of
equivalence classes $[R,\xi]$ of such 2-cocycles $(R,\xi)$ will be denoted
\[H^{2}(\KE,\underline{K}).\]

\subsection{Example.}
This example motivates the equivalence relation introduced above, and is analogous to \cite{Breen}. Let
$(\mathcal{G},\theta)$ be a gerbe on $\KE$ with band $\underline{K}$. Then $(\mathcal{G},\theta)$ gives rise to a
cocycle $(R,\xi)$ as in 11.4, which depends on the choice of objects $x_{U}$, arrows $\tilde{\sigma}$ and
homomorphisms $R(\sigma)$ for which $R(\sigma)\theta_{U}=\theta_{V}\tilde{\sigma}^{*}$.

Suppose we replace $\theta_{U}$ by $\bar{\theta}_{U}:\underline{Aut}(x_{U})\to\underline{K}(U)$. Write
$\mu_{U}=\bar{\theta}_{U}\theta_{U}^{-1}$. Then the cocycle $(R,\xi)$ is replaced by the equivalent cocycle
$(\bar{R},\bar{\xi})$ with
\[\begin{array}{rl}
\bar{R}(\sigma)&=\mu_{V}R(\sigma)\mu_{U}^{-1}\\
\bar{\xi}(\sigma,\tau)&=\mu_{W}\xi(\sigma,\tau)
\end{array}\]
(where $\sigma:V\to U$ and $\tau:U\to W$ in $\KE$ as before).

Suppose we leave $x_{U}$ and $\theta_{U}$ unchanged, but replace $\tilde{\sigma}:x_{V}\to \sigma^{*}(x_{U})$ by
$\hat{\sigma}:x_{V}\to \sigma^{*}(x_{U})$. Then $(R,\xi)$ will be replaced by the equivalent cocycle
$(\hat{R},\hat{\xi})$, defined in terms of $f(\sigma)=\theta_{V}(\hat{\sigma}^{-1}\tilde{\sigma})$, by
\[\begin{array}{rl}
\hat{R}(\sigma)&=f(\sigma)_{*}R(\sigma)\\
\hat{\xi}(\sigma,\tau)&=f(\sigma\tau)\xi(\sigma,\tau)R(\tau)(f(\sigma)^{-1})f(\tau)^{-1}.
\end{array}\]

Suppose we replace $x_{U}$ by $y_{U}$, choose arrows $\alpha_{U}:y_{U}\to x_{U}$ in $\mathcal{G}(U)$, and define
$\theta_{U}'=\theta_{U}\circ(\alpha_{U})_{*}:\underline{Aut}(y_{U})\to \underline{K}(U)$, and $\sigma':y_{V}\to
\sigma^{*}(y_{V})$ by $\sigma^{*}(\alpha_{U})\sigma'=\tilde{\sigma}\alpha_{V}$. Then also
$(\alpha_{V})_{*}(\sigma')^{*}=\tilde{\sigma}^{*}(\alpha_{U})_{*}$. One checks that the cocycle defined by these
$y_{U}$, $\sigma'$ and $\theta_{U}'$ is actually identical to the cocycle $(R,\xi)$ defined by the $x_{U}$,
$\tilde{\sigma}$ and $\theta{U}$.

This proves that the cohomology class $[R,\xi]$ associated to $\mathcal{G}$ is independent of all choices.

\subsection{Relative bands and cohomology.}
Now let $\mathcal{D}$ be a subcategory of $\KE$. A band on $(\KE,\mathcal{D})$ is a triple consisting of a band
$\underline{K}$ on $\KE$, a presheaf $\underline{K}^{\circ}$ on $\mathcal{D}$ and an isomorphism of bands
$\underline{K}^{\circ}\cong \underline{K}|\mathcal{D}$. With the obvious notion of morphisms this defines a
category of bands on $(\KE,\mathcal{D})$. An equivalent, perhaps more concrete description of this category is
obtained by taking as objects the pairs $(\underline{K},P)$ where $\underline{K}$ is a band on $\KE$ (so the
restrictions $\sigma^{*}:\underline{K}(U)\to\underline{K}(V)$ for $\sigma:V\to U$ are given as outer
homomorphisms) \textit{together} with a specific choice of an actual homomorphism
$P(\sigma):\underline{K}(U)\to\underline{K}(V)$ representing $\sigma^{*}$, for each $\sigma\in\mathcal{D}$, and
such that $P(\sigma\tau)=P(\tau)P(\sigma)$ for any composable pair $\overset{\sigma}{\gets}\overset{\tau}{\gets}$
in $\mathcal{D}$.

A \textit{relative cocycle} with coefficients in such a band $(\underline{K},P)$ is a pair $(R,\xi)$ as before,
with the additional property that $R(\sigma)=P(\sigma)$ for $\sigma\in\mathcal{D}$, and $\xi(\sigma,\tau)=1$,
whenever $\sigma$ and $\tau$ are composable arrows in $\mathcal{D}$. Equivalence is defined as before. (Note that
if $(R,\xi)$ is equivalent to $(\hat{R},\hat{\xi})$ where $\hat{R}(\sigma) =f(\sigma)_{*}R(\sigma)$ and
$\hat{\xi}(\sigma,\tau)=f(\sigma\tau)\xi(\sigma,\tau)R(\tau)(f(\sigma)^{-1})f(\tau)^{-1}$ as above, then in
particular $f$ restricts to a 1-cocycle on $\mathcal{D}$.) The set of equivalence classes of relative cocycles is
denoted
\[H^{2}(\KE,\mathcal{D};\underline{K})\]
(where we simply write $\underline{K}$ for $(\underline{K},P)$).

\subsection{Relative gerbes.}
Let $\mathcal{G}$ be a presheaf of non-empty transitive discrete groupoids on $\KE$. A \textit{section} of
$\mathcal{G}$ over $\mathcal{D}$ is a choice of objects $x_{U}\in \mathcal{G}_{U}$ (for $U\in \mathcal{D}$), which
is functorial on $\mathcal{D}$ in the sense that $x_{V}=\sigma^{*}(x_{U})$ whenever $\sigma: U\to V$ belongs to
$\mathcal{D}$. To simplify notation and anticipating the application in Theorem 12.3 below, let us assume that
$\mathcal{D}$ has the same \textit{objects} as $\KE$. Suppose now that we choose arrows
$\tilde{\sigma}:x_{V}\to\sigma^{*}(x_{U})$ in $\mathcal{G}(V)$ for each $\sigma:V\to U$ in $\KE$, with the
understanding that we choose $\tilde{\sigma}=1$ whenever $\sigma\in\mathcal{D}$. Then the groups
$\{\underline{Aut}(x_{U})\}$ form a relative band on $(\KE,\mathcal{D})$, with outer homomorphisms
$\tilde{\sigma}^{*}:\underline{Aut}(x_{U})\to\underline{Aut}(x_{V})$ defined as before and independent of
$\tilde{\sigma}$, and actual homomorphisms $\sigma^{*}$ for $\sigma\in\mathcal{D}$. We write $Band(\mathcal{G},x)$
for this relative band.

If $\theta: Band(\mathcal{G},x)\to(\underline{K},P)$ is an isomorphism of relative bands, we say
$(\mathcal{G},x,\theta)$ is a relative gerbe with band $(\underline{K},P)$. For a choice of $\tilde{\sigma}$'s as
above, we obtain a relative cocycle $(R,\xi)$ from such a relative gerbe, and hence a class
$[(\mathcal{G},x,\theta)]\in H^{2}(\KE,\mathcal{D},\underline{K})$.

\section{Extension of \'etale groupoids}
We are now ready to formulate a version of Theorem 8.2 for bundles of non-abelian Lie groups $K$.

\subsection{Extensions and bands.}
A in \S8, we consider extensions $K\to G\overset{\pi}{\to}E$ of an \'etale groupoid $E$ over $T$ by a bundle of
Lie groups $K$, the latter not necessarily abelian now. Let $\KE$ be an embedding category approximating $E$, and
let $\KE_{u}\subset\KE$ be the subcategory consisting of unit sections only (as in \S7). Suppose $\KE$ is adapted
to the extension $K\to G\to E$, so that for each arrow $\sigma:V\to U$ we can choose a lifting
$\tilde{\sigma}:V\to G$ with $\pi(\tilde{\sigma}(x))=\sigma(x)$, as before. Again, if $j:V\hookrightarrow U$
belongs to $\KE_{u}$, we agree to choose $\tilde{j}(x)=1_{x}\in G$. For each object $U\in\KE$, define
$\underline{K}(U)=\Gamma(U,K)$. Then $\underline{K}$ has the usual structure of a presheaf on $\KE_{u}$. For an
arbitrary $\sigma:V\to U$ in $\KE$, there is a map
\[\tilde{\sigma}^{*}:\underline{K}(U)\to\underline{K}(V)\]
defined by $\tilde{\sigma}^{*}(\alpha)(x)=\tilde{\sigma}(x)^{-1}\alpha(\sigma(x))\tilde{\sigma}(x)$, for
$\alpha\in\underline{K}(U)$ and $x\in V$. The corresponding outer morphism only depends on $\sigma$ and not on the
lift $\tilde{\sigma}$. Thus $\underline{K}$ has the structure of a relative band on $(\KE,\KE_{u})$ in a canonical
way, for any embedding category adapted to $K\to G\to E$. This structure replaces the action of $E$ on $K$ in the
abelian case. We say that $K\to G\to E$ is an extension with relative band $\underline{K}$.

\subsection{Non-abelian \v{C}ech cohomology.}
For an embedding category $\KE$ and a relative band $\underline{K}$ on $(\KE,\KE_{u})$, we obtain by restriction a
relative band (still denoted $\underline{K}$) on any refinement $\KE'$ of $\KE$. There is an evident restriction
map on relative cocycles, which respects equivalence, and induces a well-defined map
\[H^{2}(\KE,\KE_{u};\underline{K})\to H^{2}(\KE',\KE'_{u};\underline{K}).\]
The non-abelian \v{C}ech cohomology is now defined as the colimit over all embedding categories approximating $E$,
and written
\[\check{H}^{2}(E,T;\underline{K})=colim_{\KE}H^{2}(\KE,\KE_{u};\underline{K}),\]
analogous to the discussion in \S7. The analogue of Theorem 8.2 is

\subsection{Theorem.}
\textit{There is a 1-1 correspondence between non-abelian cohomology classes in
$\check{H}^{2}(\KE,\KE_{u};\underline{K})$ and isomorphism classes of extensions $K\to G\to E$ with band
$\underline{K}$.}
\begin{proof}
The proof is a concatenation of constructions already discussed. Indeed, there are explicit constructions of
\begin{itemize}
\item[(i)] a relative gerbe $\mathcal{G}$ on $(\KE,\KE_{u})$ with band $\underline{K}$ from an
extension $K\to G\to E$ with the same band,
\item[(ii)] a cocycle from such a relative gerbe,
\item[(iii)] an extension from such a cocycle.
\end{itemize}
We briefly discuss these constructions, and leave it to the reader to verify that these provide the required 1-1
correspondence.

(i) Consider an extension $K\to G\to E$. Let $\KE$ be an embedding category approximating $E$ and adapted to $G$.
Define a presheaf of groupoids $\mathcal{G}$ on $\KE$ as follows. Objects of $\mathcal{G}(U)$ are isomorphisms
$\rho:U\overset{\sim}{\to}W$ in $\KE$, and arrows from $\rho_{1}:U\overset{\sim}{\to}W_{1}$ to
$\rho_{2}:U\overset{\sim}{\to}W_{2}$ in $\mathcal{G}(U)$ are sections $\alpha:W_{1}\to G$ such that, for each
$x\in U$, $\alpha(\rho_{1}(x))$ is an arrow $t \rho_{1}(x)\to t\rho_{2}(x)$ in $G$ with
$\pi(\alpha(x))\rho_{1}(x)=\rho_{2}(x)$. For $\sigma:U'\to U$ in $\KE$, the restriction
$\sigma^{*}:\mathcal{G}(U)\to\mathcal{G}(U')$ is defined on objects by
\[\sigma^{*}(\rho_{1})=\rho_{1}\circ\sigma:U'\overset{\sim}{\to}W_{1}',\]
where $W_{1}'\subset W_{1}$ is the image of $\rho_{1}\circ\sigma$, and on arrows $\alpha:W_{1}\to G$ as above
simply by restricting $\alpha$ to the smaller set $W_{1}'\subset W_{1}$.

Each groupoid $\mathcal{G}(U)$ has a specific object $1_{U}:U\to U$, and $\underline{Aut}(1_{U})
\simeq\Gamma(U,K)=\underline{K}(U)$. The choice of these objects $1_{U}$ is functorial on $\KE_{u}$ but not on all
of $\KE$. Now choose liftings $\tilde{\sigma}$ for arbitrary arrows $\sigma:V\to U$ as in 12.1. These liftings
provide arrows $\tilde{\sigma}:1_{V}\to \sigma^{*}(1_{U})$ in $\mathcal{G}(V)$, conjugation by which induces a map
\[\underline{Aut}(1_{U})\overset{\sigma^{*}}{\to}\underline{Aut}(\sigma^{*}(1_{U}))\to \underline{Aut}(1_{V}),\]
sending a section $\alpha\in\underline{K}(U)$ to the section $\tilde{\sigma}^{*}(\alpha)$ as in 12.1. Thus
$(\mathcal{G},1)$ is a relative gerbe with band $\underline{K}$.

(ii) By the procedure of 11.8, this relative gerbe, together with the choice of the $\tilde{\sigma}$, defines an
explicit relative cocycle $(R,\xi)$ with $R(\sigma)=\tilde{\sigma}^{*}$ and
\[\xi(\sigma,\tau)(x)=\widetilde{\sigma\tau}(x)^{-1}\tilde{\sigma}(y)\tilde{\tau}(x)\]
where $y=t \tau(x)$, for $U\overset{\sigma}{\gets}V\overset{\tau}{\gets}W$. Notice that this is a \textit{normal}
cocycle, and that every relative cocycle can be normalized. (The discussion in 8.4 applies equally to the
non-abelian case.)

(iii) Finally, from such a normal cocycle $(R,\xi)$ we can construct an extension $K\to G\to E$ with the same
band. Indeed, the explicit construction of 8.5 again applies.
\end{proof}

\noindent
Utrecht, December 2001

\end{document}